\documentclass[12pt]{article}

\usepackage{amssymb}
\usepackage{graphicx}
\usepackage{subfigure}
\usepackage{rotating}
\usepackage{amsmath, amsthm}
\usepackage{hyperref}

\newcounter{bean}

\newcommand{\proofs}[2]{\noindent{\textbf{Proof of #1.} #2 \hfill $\square$}}
\newtheorem{theorem}{Theorem}
\newtheorem{prop} [theorem] {Proposition}
\newtheorem{lem}{Lemma}
\newtheorem{as}{Assumption}

\renewcommand{\theequation}{\arabic{section}.\arabic{equation}}

\def\cov{\mathop{\rm Cov}\nolimits}

\def\R{\mathop{\mathbb R}\nolimits}

\def\weakc{\Rightarrow}

\def\nsp{\!}

\DeclareMathOperator*{\plim}{plim} 
\input{tcilatex}

\begin{document}

\begin{titlepage}
\center{\LARGE\mbox{Specification~Tests~for~Nonlinear~Dynamic~Models}}
\center{\textsc{Igor Kheifets}\\New Economic School, Moscow\\ 
October 13, 2014}
\begin{abstract}
We propose a new adequacy test and a graphical evaluation tool for nonlinear  dynamic models. The proposed techniques can be applied in any setup where  parametric conditional distribution of the data is specified, in particular to models involving conditional volatility, conditional higher moments, conditional quantiles, asymmetry, Value at Risk models, duration models, diffusion models, etc.  Compared to other tests, the new test properly controls the nonlinear dynamic behavior in conditional distribution and does not rely on smoothing techniques which require a choice of several tuning parameters. The test is based on a new kind of multivariate empirical process of contemporaneous and lagged probability integral transforms. We establish weak convergence of the process under parameter uncertainty and  local alternatives.  We justify a  parametric bootstrap approximation that accounts for parameter estimation effects often ignored in practice.  Monte Carlo experiments show that the test has good finite-sample size and power properties. Using the new test and graphical tools we check the adequacy of various popular heteroscedastic models for stock exchange index data.
\end{abstract}
\begin{flushleft}
\end{flushleft}
\thispagestyle{plain}
\end{titlepage}

\setcounter{page}{2}
\section{INTRODUCTION}
\setcounter{theorem}{0}\setcounter{equation}{0}

In this paper we develop a methodology for testing the goodness-of-fit of a parametric conditional distribution in nonlinear time series model.  More precisely,
let  $Y_{t},\ t=0, \pm 1,\pm 2,\ldots$ be a univariate real-valued time series. Let $\Omega _{t}$ be a sigma-field generated by the observations obtained up to time $t$
(i.e. by $Y_{t-1},Y_{t-2},\ldots$, information set at time $t$, not including $Y_{t}$).
We consider the family of
conditional distribution functions $F_{t}\left(\cdot\mid\Omega _{t},\theta \right)$, parameterized by $\theta
\in \Theta $, where $\Theta \subset \R^{L}$ is a finite dimensional
parameter space.
We permit changes over time in the
functional form of the distribution using subscript $t$ in $F_{t}$. Sometimes, with a small abuse of notation,  we will use   $F_{t}\left(\cdot\mid\theta\right)$ for $F_{t}\left(\cdot\mid\Omega _{t},\theta\right)$.
We aim to test the correct specification of the conditional distribution against the general alternative. Our null hypothesis is

$H_{0}$ : The conditional distribution of $Y_{t}$ conditional on $\Omega
_{t} $ is in the parametric family $F_{t}\left(\cdot\mid\Omega _{t},\theta \right)$ for some $%
\theta _{0}\in \Theta $.

Testing the specification of nonlinear time series models is crucial in applied statistics,
macroeconomics and finance for making relevant analysis. It is often
not enough to check only conditional moments.  It is often necessary to check the specification of the 
\textit{conditional distribution function}, which is equivalent to jointly specifying all conditional
features of a process, including the conditional mean, variance, and quantiles. 
For instance,
knowing the true conditional distribution is essential for applying \textit{%
efficient} maximum likelihood (ML) methods to many models. The importance of nongaussian likelihood methods is stressed in 
 Harvey (2010) and Qi, Xiu and Fan (2010). 
The conditional distribution is linked to the hazard function, e.g., in the autoregressive conditional duration (ACD) model of Engle and Russel (1998). The knowledge of the conditional
distribution provides Value at Risk (VaR),  defined as a quantile of the return distribution, and is used to assess risk in finance. 
The description of the return distribution in its whole support is especially important in risk management to estimate the downside risk of nonlinear portfolios.
Tests for   conditional distribution may serve  to evaluate the density forecasts of macroeconomic variables such as inflation and risks in finance and
insurance, see Diebold et al. (1998), Thomson (2008).

The class of models which specify parametric conditional distribution is very broad. In dynamic location-scale models, such as ARMA and GARCH,
conditional distribution is simplified to an unconditional distribution
function of innovations, typically normal or student-$t$, and a dependence structure, which is modeled
only through the conditional first two moments.  If $Y_{t}=\mu
_{t}(\alpha )+\sigma _{t}(\beta )\varepsilon _{t},$ where $\mu _{t}(\alpha )$
and $\sigma _{t}(\beta )$ are measurable with respect to $\Omega _{t}$ and $%
\varepsilon _{t}$ are independent identically distributed (iid) with cumulative distribution function (cdf) $F_{\varepsilon t}(\varepsilon ;\gamma
)$, the conditional distribution can be expressed  as 
$F_{t}\left(y|\Omega _{t},\theta\right)=F_{\varepsilon t}\left(\left({y-\mu_t(\alpha)}\right)/{\sigma_t(\beta)};\gamma\right)$,
with $\theta=\left(\alpha,\beta,\gamma\right)$.
Another example of dynamic location-scale models is a data series
discretely sampled from the Ornstein-Uhlenbeck process, which is commonly used to
model interest rates. Usually, these models are examined by testing dynamics and marginal distributions of residuals separately, but these tests are inconsistent against many alternatives to $H_0$. Furthermore, for nonlinear  models it is often difficult to obtain residuals, while the  conditional distribution of the observations is easily specified.
There are 
examples beyond the class of dynamic location-scale
models. Discretely sampled series from a diffusion model (e.g. from the mean-reverting
square-root process, Feller process) may have a closed form
conditional distribution, which does not necessary belong to the
location-scale class. %
Even if sampled data from a diffusion model does not
have conditional distribution, it may be approximated.
The stochastic volatility (SV) models have no explicit form of conditional
distribution, but it can be simulated in a manner similar to SV simulated maximum
likelihood estimation.

Our work relies on the following well known  result for the (conditional) \emph{probability integral transform} (also called \emph{quantile transform}), which is similar in spirit to Rosenblatt (1952).

\begin{prop}
\label{propFY} Suppose that the conditional distributions $F_{t}(y|\Omega _{t},\theta )$ are
 continuous and
strictly increasing in $y$. If $Y_t\mid\Omega _{t}\sim F_{t}\left(y|\Omega _{t},\theta_0\right)$, then random variables 
 $U_{t}=F_{t}\left(Y_{t}|\Omega _{t},\theta_0\right)$ are iid standard uniform.
\end{prop}

Random variables $U_{t}$, sometimes called 
\emph{generalized residuals}\footnote{E.g. Hong and Li (2005), p. 39.}, are easy to obtain.
In the location-scale model, the transform is
$U_{t}=F_{\varepsilon t}\left({\left(Y_{t}-\mu_t(\alpha_0)\right)}/{\sigma_t(\beta_0)};\gamma_0\right)$,
$\theta_{0}=\left(\alpha_0,\beta_0,\gamma_0\right)$.
In continuous time models, i.e. when $Y_{t}$ comes from $dY_{t}=\mu_{t}(Y_{t},\theta)dt+\sigma _{t}(Y_{t},\theta)dW_{t}$, where $%
\mu_{t}(Y_{t},\theta)$ and $\sigma _{t}(Y_{t},\theta)$ are the drift and diffusion
functions respectively, and $W_{t}$ is standard Brownian motion, we can
test for the correct transition density family $\{p(x,t|y,s;\theta )\}$,
the conditional density of $Y_{t}=x$ given $Y_{s}=y,$ $s<t$. For any $\Delta>0$ the transform is defined as $U_{t}=\int_{-\infty }^{Y_{t\Delta
}}p(x,t\Delta |Y_{(t-1)\Delta },(t-1)\Delta ;\theta )dx$. 

Similarly to residuals of a linear regression model, working with generalized residuals is easier than with a conditional distribution of arbitrary form, and researchers often apply existing statistical tests to the marginal distribution and/or dynamic structure of generalized residuals of nonlinear models, see e.g. Kalliovirta (2012) and references therein. For example, 
the uniformity of $U_t$ is checked by transforming them into Gaussian random variables (under the null) and using the normality test of Jarque and Bera (1987) or the residual empirical process (Bai 2003), QQ plots (Diebold et al. 1998) and their dynamics through correlations of residuals and squared residuals (Box and Pierce, 1970; Ljung and
Box, 1978; Diebold et al., 1998; Du and Escanciano, 2013) 
or serial independence test (Skaug and Tj\o stheim 1993). Although these procedures can give insight into the model misspecification, they check only specific properties of a model and are inconsistent for $H_0$.  Typically they are applied sequentially and size is not corrected for the sequential testing problem. Moreover, sometimes existing methods are not valid.  A typical mistake is to apply Kolmogorov --- Smirnov type tests to verify iid-ness and functional form of marginal distribution, while this test is designed for verifying the marginal distribution of iid random variables and should not be used with generalized residuals. Since $H_0$ requires examining simultaneously the uniformity \emph{and}
independence of the generalized residuals, not uniformity under independence, the Kolmogorov-Smirnov test for $U_t$,  or test of Bai (2003),  do not control the dynamics in $U_t$ and miss important deviations from the null, which is shown analytically and in simulations below. 

In this paper we propose a test which is specifically designed for $H_0$.
Our test is based on a new multivariate empirical process which incorporates the difference between the empirical
joint distribution of lags of $U_t$ and the product of the uniform marginals. The asymptotic properties of such a process have never been studied before;  therefore, we establish weak convergence results for the underlying process that allow us,
under standard conditions,
to prove consistency and derive asymptotic properties, including the
distribution under root-$n$ local alternatives, of such tests taking into
account parameter estimation effect. Since the asymptotic distribution is case
dependent, critical values can not be tabulated, and we prove that the
parametric bootstrap distribution approximation is valid. Bootstrap implementation is straightforward  since the DGP is known under $H_0$. Since   no smoothing
techniques are used,  there are no user chosen parameters, and
we have standard root-$n$ rates of convergence.  The underlying empirical process is piece-wise quadratic; therefore, computation of the test statistics does not hinge on numerical optimization or numerical integration.  The only model-depended routines needed are simulation and estimation of the original model.

The Monte Carlo study shows that we have good size and power
properties for both linear and nonlinear models. 
We show both analytically and with simulations that our test achieves better performance compared to tests based on a univariate residual empirical process as in Bai (2003). An explicit
closed form expression for our test statistics that avoids numeric approximation is available.
We also propose a graphical tool of model evaluation in order to identify  potential sources of misspecification.

The importance of  simultaneous checks for the uniformity and
independence of $U_t$ has been emphasized in series of papers (Berkowitz 2001, Gonzalez-Rivera et al. 2011, Kalliovirta 2012); however, they test correlations instead of independence of $U_t$ and  are  inconsistent against a higher order dependence. Hong and Li (2005) examines uniformity and
independence by comparing the unrestricted kernel estimator of the joint density of $U_t$ with the model implied one. In contrast to our test, their test has a slower rate of convergence, cannot distinguish Pittman local alternatives and depends on bandwidth selection procedures, which may be impractical. Our test is close in spirit to the empirical process based goodness-of-fit tests of $H_0$ by Andrews (1997), Rothe and Wied (2013) and Delgado and Stute (2008), although their results hinge on the iid-nes of the original data, while our test not only can be applied to time series, but also is able to detect misspecification in dynamic structure.

The rest of the paper is organized as follows.  The new test is introduced in Section~2. Asymptotic properties of the test 
and bootstrap justification are provided in Section~3.
Monte Carlo experiments and application to the daily stock index are in Section~4. In Section~5, we briefly discuss how results in this paper can be used for testing the specification of multivariate models, and a conclusion is offered in Section~6. Proofs  are deferred to  the Appendix.

\section{THE NEW TESTS}
\setcounter{theorem}{0}\setcounter{equation}{0}
Our goal is to measure how far the generalized residuals are from being independent and uniform. For instance, under the null  for  $r=(r_{1},r_{2})\in {[0,1]}%
^{2}, $
\begin{equation}\label{eq:p12}
P\left({U}%
_{t}\leq r_{1},{U}_{t-1}\leq r_{2}\right)=r_{1}r_{2}.
\end{equation} 
 Then the nonparamtric estimator of the {\em joint} distribution of $(U_t,U_{t-1})$ must be close to the bivariate uniform. This motivates us to
consider the following empirical process
\begin{equation*}
{V}_{2n}(r)=\frac{1}{\sqrt{n-1}}\sum_{t=2}^{n}\left[ I({U}_{t}\leq r_{1})I({U%
}_{t-1}\leq r_{2})-r_{1}r_{2}\right], \label{eqV2n}
\end{equation*}%
where $I(\cdot)$ is the indicator function.
This bi-parameter empirical
process  incorporates the difference between the empirical
joint distribution and the product of the uniform marginal (cumulative) distribution functions of contemporaneous
 and lagged generalized residuals, i.e. examines 
{simultaneously} the uniformity and
independence of $U_t$. Note, the process $V_{2n}$ checks the 
{implication} of the $H_0$ stated in Equation (\ref{eq:p12}). 
Empirical process theory is useful in the context of goodness-of-fit testing.  For example, Bai (2003) used a univariate empirical process to test $H_0$; however, his test is inconsistent (see Corradi and Swanson (2006) and our discussion below). In an iid setup, Delgado and Stute (2008) solve the issue with Bai's test using a similar bivariate process; however, their results do not hold in the dynamic setting we consider in this paper.

For an illustration, let $Y_t$ be negative stock returns (losses) and consider Value at Risk,  defined as a quantile of the (conditional) distribution function, $VaR_t(r_1)=F^{-1}\left(r_1|\Omega_t,\theta_0 \right)$. Our process backtests Value at Risk, i.e.,  it verifies that events $A_i=\{Y_{t+1-i}\le VaR_{t+1-i}(r_i)\}$, $i=1,2$ (losses do not exceed Value at Risk),  are independent and have the right (unconditional) coverage, for all  $(r_{1},r_{2})\in {[0,1]}%
^{2}$, thus ensuring \emph{conditional} coverage. 
For the recent literature on Value at Risk backtesting see Escanciano and Olmo (2011) and references therein, where the danger of using only unconditional checks  is emphasized (as in Bai, 2003).

In practice, we rarely know $\theta _{0}$ either $\{Y_{t},t\leq 0\}$. We can
approximate $U_{t}$ with $\hat{U}_{t}=F_{t}(Y_{t}|\tilde{\Omega}_{t},{\hat{%
\theta}})$ where ${\hat{\theta}}$ is an estimator of $\theta _{0}$ and the
truncated information $\tilde{\Omega}_{t}$ is based on $\{Y_{t-1},Y_{t-2},\ldots ,Y_{1}\}$ and write
\begin{equation}
\hat{V}_{2n}(r)=\frac{1}{\sqrt{n-1}}\sum_{t=2}^{n}\left[ I(\hat{U}_{t}\leq
r_{1})I(\hat{U}_{t-1}\leq r_{2})-r_{1}r_{2}\right] .  \label{eqVhat2n}
\end{equation}
The process $\hat{V}_{2n}(r)$ measures the distance to the null hypothesis
for each $r$, so we need to choose a metric in $[0,1]^{2}$ to aggregate
for all $r$. For any continuous functional $\Gamma (\cdot )$ from the set of uniformly bounded real functions on $%
[0,1]^{2}$, to~$\mathbb{R}$,
\begin{equation*}
D_{2n}=\Gamma (\hat{V}_{2n}(r)).
\end{equation*}%
In particular we consider Cramer-von Mises (CvM) and Kolmogorov-Smirnov (KS) statistics
\begin{equation*}
D_{2n}^{CvM}=\int_{[0,1]^{2}}\hat{V}_{2n}(r)^{2}dr\text{ and }%
D_{2n}^{KS}=\sup_{[0,1]^{2}}\left\vert \hat{V}_{2n}(r)\right\vert .
\end{equation*}
Note that numerical integration and maximization can be avoided since the process $\hat{V}_{2n}(r)$ is piece-wise quadratic in $r$ .

Recall that (\ref{eqVhat2n}) checks only pairwise dependence. To check $p$%
-wise independence (see Delgado, 1996) we can write
\begin{equation*}
{V}_{pn}(r)=\frac{1}{\sqrt{n-p+1}}\sum_{t=p}^{n}\left[
\prod_{j=1}^{p}I({U}_{t-j+1}\leq r_{j})-r_{1}r_{2}\ldots r_{p}\right]
\end{equation*}
and use the test statistics
\begin{equation*}
D_{pn}^{CvM}=\int_{[0,1]^{p}}\hat{V}_{pn}(r)^{2}dr\text{ and }%
D_{pn}^{KS}=\max_{[0,1]^{p}}\left\vert \hat{V}_{pn}(r)\right\vert.
\label{eq:Dpn}
\end{equation*}
For instance, the process $\hat{V}_{1n}$ compares the nonparametric estimator of the {\em marginal} distribution of $U_t$ with the uniform, formalizing QQ-plots in Deibold et al. (1998) (see Bai (2003) for a related test). Interestingly, the following decomposition holds
\begin{equation*}
V_{2n}\left( r\right)  =\sqrt\frac{{n}}{{1-n}}\left\{V_{1n}\left( r_{1}\right) r_{2}+V_{1n}\left(
r_{2}\right) r_{1}\right\} +V_{2n}^C\left( r\right),
\end{equation*}
where
\begin{equation*}
V_{2n}^C\left( r\right)=\frac{1}{\sqrt{n-1}}\sum_{t=2}^{n}\left[ I({U}_{t}\leq r_{1})-r_{1}\right]
\left[ I({U}_{t-1}\leq r_{2})-r_{2}\right].
\end{equation*}
Thus, we have a combination of a   univariate empirical process $V_{1n}$ which
accounts for the unconditional/marginal distribution fit,  and the correlation of the
indicators, for dynamic check. Note also, that $V_{2n}(r_1,1)=\sqrt{n/(n-1)}V_{1n}(r_1)$.

To test $j$-lag ($j=1,2,...$) pairwise independence, define process
\begin{equation*}
\hat{V}_{2n,j}(r)=\frac{1}{\sqrt{n-j}}\sum_{t=j+1}^{n}\left[ I(\hat{U}%
_{t}\leq r_{1})I(\hat{U}_{t\text{-}j}\leq r_{2})-r_{1}r_{2}\right],
\label{eq:V2nj}
\end{equation*}%
and test statistics 
\begin{equation}
D^{CvM}_{2n,j}=\int_{[0,1]^{2}}\hat{V}_{2n,j}(r)^{2}dr\text{ and }%
D^{KS}_{2n,j}=\max_{[0,1]^{2}}\left\vert \hat{V}_{2n,j}(r)\right\vert.
\label{eq:Dnj}
\end{equation}%
Likewise, the Portmanteau test (e.g. Ljung and Box 1978) uses information on autocorrelations in a number of lags. Therefore, 
we can aggregate (\ref{eq:Dnj}) across $j$ using the sum or maximum operator, so that for $k<n$, we get test statistics
\begin{equation*}
ADJ_{kn}=%
\sum_{j=1}^{k}D_{2n,j}^{CvM}\text{ and }MDJ_{kn}=%
\max_{j=1,...,k}D_{2n,j}^{KS}.  \label{eq:ADn}
\end{equation*}
To use the information from statistics $D_{1n}$, which provides additional explicit  account for the unconditional distribution misspecification, we introduce 
\begin{equation*}
ADJ_{kn}^0=%
D_{1n}^{CvM}+ADJ_{kn}\text{ and }
MDJ_{kn}^0=%
\max\left(D_{1n}^{KS},MDJ_{kn}\right).  \label{eq:ADn0}
\end{equation*}
Many other combinations of the
statistics $D_{2n,j}$ and $D_{pn}$,
including aggregation across $p$ or summing with different weights, can be considered.
The optimal choice of aggregation, in particular of $k,p$ and the weighting scheme, possibly data-driven, will depend on the alternative and is not addressed here. We note that in practice one should add lags with care since including too many lags may reduce power in small samples. As we show below with simulations, pairwise statistics with one or five equal-weighted lags  already captures many relevant alternatives in moderate sample sizes. 

We propose the following  interpretation of the introduced statistics.
In a nonlinear/nongaussian setup, correlations and Ljung-Box statistics
do not provide the whole picture of the goodness-of-fit, since  dynamics may not result in significant
correlations. Therefore we propose to use 
$D_{2n,j}$ as \emph{generalized sample autocorrelations}
and $ADJ_{kn}^0$ or $MDJ_{kn}^0$  as \emph{generalized} Ljung-Box tests. Thus we can consider a
\emph{generalized autocorrelogram}, drawing $D_{2n,j}$ against $j$.
The difference is that there may be dependence across different $j$ (which has to be taken into account
if one wants to make joint inference) and that all values
are nonnegative.\footnote{Marginal test statistics, such as Jarque-Bera or Kolmogorov-Smirnov and  autocorrelations/Ljung-Box often are used together in specification testing. 
Unlike their traditional counterparts, generalized sample autocorrelations and generalized Ljung-Box test statistics introduced here incorporate information from marginals test statistics (by setting $r_2=1$). Therefore our procedure requires no additional marginal check, avoiding sequential testing. 
}
In our Monte Carlo study, we compare the performance of traditional and generalized statistics.

As we will see in the next section (Proposition~\ref{proplimitV}), when no parameters are estimated, our test statistics are asymptotically distribution-free, and their critical values can be simulated and tabulated. In practice, however, the parameter estimation effect needs to be taken into account, which requires a bootstrap method to approximate the distribution of the test statistics.
Under $H_{0}$ we know the parametric conditional
distribution; therefore, we apply a parametric bootstrap to mimic the $H_{0}$
distribution based on $F_{t}(\cdot |\cdot,\hat{%
\theta})$, which is essentially the same as Monte Carlos simulations. We introduce the algorithm now.

\setcounter{bean}{0}
\begin{center}
  \begin{list}
{\textsc{Step} \arabic{bean} .}{\usecounter{bean}}
\item Estimate the model with the original data $Y_{t}$, $t=1,2,...,n$, get parameter
estimator $\hat{\theta}$, and get test statistic $\Gamma (\hat{V}_{2n})$.

\item Simulate $Y_{t}^{\ast }$ with $F_{t}(\cdot |\Omega _{t}^{\ast },\hat{%
\theta})$ recursively for $t=1,2,...,n$, where $\Omega _{t}^{\ast
}=(Y_{t-1}^{\ast },Y_{t-2}^{\ast },...)$.\label{basimulate}

\item Estimate the model with simulated data $Y_{t}^{\ast }$, and get $\theta ^{\ast
}$, get bootstrapped statistics $\Gamma (\hat{V}_{2n}^{\ast })$.

\item Repeat 2-3 $B$ times, and compute the percentiles of the empirical
distribution of the $B$ bootstrapped statistics.

\item Reject $H_{0}$ if $\Gamma (\hat{V}_{2n})$ is greater than the $(1-\alpha )$%
th percentile of the empirical distribution.
\end{list}
\end{center}

This technique is much easier to implement than block bootstrap techniques often utilized in a time series setup (e.g. Corradi and Swanson  2006) because it does not require a block length choice.  Note that the block bootstrap test of Corradi and Swanson (2006) cannot be directly applied here, since they test a different null hypothesis which allows dynamic misspecification. Their null hypothesis does not necessarily provide iid transforms, and thus they do not control dependence in transforms, which is crucial for the power of our test. One case when the block bootstrap is desirable is when model simulation is very costly computationally.  Indeed, while the block bootstrap method requires an estimation of a model as many times as the parametric bootstrap, generating bootstrapped samples from blocks may be faster than simulation from the distribution for some nonlinear models. However, usually estimation is more costly that simulation. In particular, for models considered in our Application section, typical estimation takes 50 longer than simulation which takes 0.0164 seconds, thus possible speed improvement from the block bootstrap method is very small. In reality, we need also to compute test statistics for each bootstrap sample; therefore, the relative simulation cost is even less.\footnote{
Methods have been proposed in the literature to make empirical process based tests distribution-free even in the presence of the parameter estimation effect. 
For instance, to make a distribution-free specification test, Delgado and Stute (2008) use the Khmaladze transform for their bivariate process. This method has its own limitations: it requires analytical derivation, estimation  and programming of the transform for each model, but it may be useful for models which are computationally hard to estimate and where bootstrap methods may take a very long time. 
Because of time dependence, the method of Delgado and Stute (2008) cannot be applied to our case. The extension of the Khmaladze transform to our case is left for future research.}

\section{ASYMPTOTIC PROPERTIES}
\setcounter{theorem}{0}\setcounter{equation}{0}\setcounter{as}{0}

In this section, we derive the asymptotic properties of the proposed statistics. When $p=1$, it can be done using standard arguments for univariate empirical processes. However for $p>1$, the theory is substantially different. The difficulty is that lags of the same variables enter in different dimensions of the processes. Here we discuss in detail the case $p=2$, and results for $p=1$ follow by simply fixing $r_2=1$. Generalization for $p>2$ can be done along the same lines, however it is lengthly and thus omitted. Note that the case $p=2$ is the most relevant for practical purposes.

We start
with a simple case in which parameters are known, then study how  the
asymptotic distribution changes if  parameters are estimated. We provide analysis under
the null, under local and fixed alternatives. We impose assumptions
on conditional cdf, the form of parametric family of cdf, dynamics and on the
estimator. 

\begin{as} 
{The conditional distributions $F_{t}(y|\Omega _{t},\theta )$ are
 continuous and
strictly increasing in $y$ for all $\theta\in\Theta$.}
\end{as}

We first describe the asymptotic behavior of the process $V_{2n}$ under
$H_0$. Denote by ``$\weakc$" weak convergence of stochastic processes as random elements of
the Skorokhod space $D\left([0,1]^2\right)$ and by ``$a \wedge b$" a minimum between $a$ and $b$.

\begin{prop}
\label{proplimitV} Suppose Assumption 3.1 holds. Then under $H_0$
\begin{equation*}
V_{2n}\weakc V_{2\infty },
\end{equation*}%
where $V_{2\infty }({r})$ is bi-parameter zero mean Gaussian process with
covariance
\begin{equation}
\mathop{\rm Cov}\nolimits_{V_{2\infty }}({r},{s})=(r_{1}\wedge
s_{1})(r_{2}\wedge s_{2})+(r_{1}\wedge s_{2})r_{2}s_{1}+(r_{2}\wedge
s_{1})r_{1}s_{2}-3r_{1}r_{2}s_{1}s_{2}.  \label{eqasscov}
\end{equation}
\end{prop}
The covariance is different from that of the two-parameter Brownian Bridge
\begin{equation*}
\lim_{n\rightarrow \infty }\mathop{\rm Cov}\nolimits_{S_{n}}({r},{s}%
)=((r_{1}\wedge s_{1})-r_{1}s_{1})((r_{2}\wedge s_{2})-r_{2}s_{2}),
\end{equation*}%
the asymptotic distribution from Skaug and
Tj\o stheim (1993). Note that if we fix $%
r_{2}=s_{2}=1,$ we get one parameter Brownian Bridge, which establishes the limit of $V_{1n}$. %

Suppose the
conditional distribution function $H_{t}(y|\Omega _{t})$ is not in the
parametric family $F_{t}(y|\Omega _{t},\theta )$, i.e. for each $\theta \in
\Theta $ there exists 
$y\in \mathop{\mathbb R}\nolimits$ and $ {t_{0}}$ so that it occurs
with positive probability 
$H_{t_{0}}(y|\Omega _{t_{0}})\neq F_{t_{0}}(y|\Omega _{t_{0}},\theta )$. For
any $0<\delta /\sqrt{n}<1$ define the conditional cdf
\begin{equation*}
G_{nt}(y|\Omega _{t},\theta )=\left( 1-\frac{\delta }{\sqrt{n}}\right)
F_{t}(y|\Omega _{t},\theta )+\frac{\delta }{\sqrt{n}}H_{t}(y|\Omega _{t}).
\end{equation*}%
and the local alternatives
\bigskip

$H_{1n}$: The conditional distribution function of $Y_{t}$ is equal to $G_{nt}(y|\Omega _{t},\theta
_{0})$.
\bigskip

We extend weak convergence in Proposition \ref{proplimitV} under $H_{1n}$ and  
to the case of the composite null hypothesis,
for which we study the parameter estimation effect (Durbin, 1973):  how $%
\hat{V}_{2n}({r})$ differs from $V_{2n}({r})$. Let $\Vert \cdot \Vert $ denote
Euclidean norm for matrices: $\Vert A\Vert =\sqrt{\mathop{\rm tr}%
\nolimits(A^{\prime }A)}$ and for $\varepsilon >0,$ $B(a,\varepsilon )$\ is
an open ball in $R^{L}$ with the center in the point $a$ and the radius $%
\varepsilon $. In particular, for some $M>0$ denote $B_{n}=B\left( \theta
_{0},Mn^{-1/2}\right) =\{\theta :||\theta -\theta _{0}||\leq Mn^{-1/2}\}$.
For simplicity, in the following  we assume no information truncation, $\tilde\Omega_t=\Omega_t$, omit $\Omega_t$ and write $\eta_t\left(r,u,v\right)=F_{t}\left( F_{t}^{-1}\left(r |u\right)
|v\right)$.  This requirement can be relaxed, e.g., when Condition A4 of Bai (2003) is satisfied, the difference between empirical processes $\hat V_{2n}$ with and without information truncation is $o_p(1)$, and thus our results are still valid  for GARCH process. We will need the following assumptions.

\begin{as}
\begin{itemize}
\item[(a)] Suppose that
\begin{equation*}
E\sup_{t=1,..,n}\sup_{u\in B_n}\sup_{r\in [0,1]}\left\vert \eta _{t}\left( r,u,\theta_0\right)-r \right\vert
=O\left( n^{-1/2}\right).
\end{equation*}
\item[(b)]  $\forall M\in(0,\infty)$, $\forall M_2\in(0,\infty)$ and $\forall\delta> 0 $
\begin{equation*}
\sup_{r\in [0,1]}
\frac{1}{\sqrt{n}}
\sum_{t=1}^{n}
\sup_{
\substack{
||u-v||\le
 M_2 n^{-1/2-\delta}\\
u,v\in B_n
}
 }
\left|
{\eta
}_{t}\left( r,u,\theta_0\right)
-
{\eta }_{t}\left(
r,v,\theta_0\right)
\right| 
=o_{p}\left(1\right).
\end{equation*}
\item[(c)]  $\forall M\in(0,\infty)$, $\forall M_2\in(0,\infty)$ and $\forall\delta> 0 $
\begin{equation*}
\sup_{|r-s|\le M_2 n^{-1/2-\delta} }
\frac{1}{\sqrt{n}}
\sum_{t=1}^{n}
\sup_{u\in B_n}
\left|
{\eta
}_{t}\left(r, u,\theta_0\right)
-
{\eta }_{t}\left(s, u,\theta_0\right)
\right| 
=o_{p}\left(1\right).
\end{equation*}
\item[(d)]  $\forall M\in(0,\infty)$, 
there exists a uniformly continuous (vector) function $h(r)$
from $[0,1]^{2}$ to $R^{L}$, such that
\begin{equation*}
\sup_{u\in B_n }\sup_{r\in \lbrack
0,1]^{2}}\left\vert
\frac{1}{\sqrt{n}}\sum_{t=2}^{n}h_t-h(r)^{\prime }{\sqrt{n}\left( u-\theta_0\right) }%
\right\vert =o_{p}(1).
\end{equation*}
where
\begin{equation*}
h_t=
\left(\eta_{t-1}\left( r_{2},u,\theta_0\right)
-r_{2}\right) r_{1} 
+\left( \eta_{t}\left( r_{1},u,\theta_0\right) -r_{1}\right)
I\left( F_{t-1}\left(  Y_{t-1}|u\right)\leq r_{2}\right).
\end{equation*}
\end{itemize}
\end{as}
We impose two types of restrictions. 
The first is the smoothness of the distributions with respect to parameters and data.  Similar assumptions have been used previously in statistical literature (e.g., Loynes, 1980). 
The second is the convergence in probability of certain averages and  implicitly imposes restrictions on the data dynamics and can be established by means of ULLN. 
This part can be also verified directly given a particular model, see e.g. proof of Theorem 3 in Bai (2003), but in general it can be a difficult task. 
If the cdf $F_{t}\left( x|\theta \right) $ is continuous differentiable with
respect to $\theta $ uniformly in $t$ and $x$, then by the mean value
theorem there exists $v^{\star }$ on the segment between $u$ and $v$, possibly
depending on $t,$ $r,$ such that $F_{t}\left( F_{t}^{-1}\left(
r_{1}|u\right) |v\right) -r_{1}=\nabla _{\theta }F_{t}\left(
F_{t}^{-1}\left( r_{1}|u\right) |v^{\star }\right) {\left( u-v\right) .}$
Therefore the following Conditions (a') and (b') (which are standard, see Bai, 2003)  
are sufficient for (a)-(c).
\textit{\begin{itemize}
\item[(a')]  
$\forall M\in(0,\infty)$, there exists a uniformly (in $x$ and $t$) continuous (with
respect to $\theta $) gradient $\nabla _{\theta }F_{t}\left( x|\theta
\right) $ which is also uniformly bounded:
\begin{equation*}
\sup_{\theta \in B_n}\sup_{x,t}\left\Vert
\nabla _{\theta }F_{t}\left( x|\theta \right) \right\Vert =O_{p}\left(
1\right) .
\end{equation*}
\item[(b')]  
$\forall M\in(0,\infty)$, there exists a density $f_{t}\left( x|\theta
\right) $ which is also uniformly bounded:
\begin{equation*}
\sup_{\theta \in B_n}\sup_{x,t}
f_{t}\left( x|\theta \right) =O_{p}\left(
1\right) .
\end{equation*}
\end{itemize}
Condition (d) holds if we add the following assumption (c') with
 $h\left( r\right) =h_{1}\left( r_{2}\right) r_{1}+h_{2}\left(
r\right).$ 
\begin{itemize}
\item[(c')] 
$\forall M\in(0,\infty)$, 
there exist uniformly continuous (vector) functions $h_1(r)$ and  $h_2(r)$
from $[0,1]^{2}$ to $R^{L}$, such that
\begin{equation*}
\sup_{u,v\in B_n}\sup_{r\in \lbrack
0,1]}\left\Vert \frac{1}{n}\sum_{t=1}^{n}
 \nabla _{v}\eta_{t}\left( r,u,v\right) 
-h_1(r)\right\Vert =o_{p}(1)
\end{equation*}
and
\begin{equation*}
\sup_{u,v\in B_n }\sup_{r\in \lbrack
0,1]^{2}}\left\Vert
\frac{1}{{n}}\sum_{t=2}^{n}h_{2t}-h_2(r)
\right\Vert =o_{p}(1),
\end{equation*}
where
\begin{equation*}
h_{2t}=
 \nabla _{v}\eta_{t}\left( r_{1},u,v\right) 
\left\{I\left( F_{t-1}\left(  Y_{t-1}|u\right)\leq r_{2}\right)-r_2\right\}.
\end{equation*}
\end{itemize}
}
Conditions similar to (c') have been
used in the empirical process literature (Bai, 2003).  For $p=1$ (and also in the iid case and any $p$), it is enough to use
only the condition for $h_1$. For $p=2$, we need
an additional condition for $h_2$. In the iid case, (c') holds automatically with $h_1(r)=\nabla _{v}\eta_{1}\left( r,u,v\right) 
$.  In a stationary and ergodic case, (c') holds by taking the unconditional expectation and applying ULLN. In a dynamic heterogeneous case, ULLN also exist (Potscher and Prucha, 1997).

The term $h\left( r\right) ^{\prime }${$\sqrt{n}%
\left( \hat{\theta}-\theta _{0}\right) $ will appear in the expansion of }$%
\hat{V}_{2n}({r})$ around $V_{2n}({r})$ and will reflect the parameter
estimation effect. 
Hence, to identify the limit of $\hat{V}_{2n}(r)$ we need to study the
limiting distribution of random vector $\left( {V_{2n}(r),\sqrt{n}}\left( {%
\hat{\theta}-\theta _{0}}\right) ^{\prime }\right) ^{\prime }$.  We make assumptions
on the estimation procedure.

\begin{as}
Under $H_{1n}$, the estimator $\hat{\theta}$ admits a linear expansion
\begin{equation}
\sqrt{n}(\hat{\theta}-\theta _{0})=\delta\mu+\frac{1}{\sqrt{n}}\sum_{t=1}^{n}\psi\left(Y_t,
\Omega _{t},\theta_0\right)+o_{p}(1),  \label{eqestlinear}
\end{equation}%
with $E\left( \psi\left(Y_t,\Omega _{t},\theta_0\right)|\Omega _{t}\right) =0$ and
\begin{equation*}
\frac{1}{n}\sum_{t=1}^{n}\psi\left(Y_t,\Omega _{t},\theta_0\right)\psi\left(Y_t,\Omega _{t},\theta_0\right)^{\prime
}
\stackrel{p}{\rightarrow } \Psi .
\end{equation*}
\end{as}
This assumption is satisfied for ML and nonlinear least square (NLS) estimators under minor
additional conditions. It will allow to apply the CLT for random vector  $({%
V_{2n}(r),}\frac{1}{\sqrt{n}}\sum_{t=1}^{n}\psi\left(Y_t,\Omega _{t},\theta_0\right)^{\prime })^{\prime }$.
Define
\begin{equation*}
C_{n}(r,s,\theta_0 )=E\left(
\begin{array}{c}
{V_{2n}(r)} \\
\frac{1}{\sqrt{n}}\sum_{t=1}^{n}\psi\left(Y_t,\Omega _{t},\theta_0\right)%
\end{array}%
\right) \left(
\begin{array}{c}
{V_{2n}(s)} \\
\frac{1}{\sqrt{n}}\sum_{t=1}^{n}\psi\left(Y_t,\Omega _{t},\theta_0\right)%
\end{array}%
\right) ^{\prime }
\end{equation*}%
and let $(V_{2\infty }(r),\psi _{\infty }^{\prime })^{\prime }$ be a zero
mean Gaussian process with covariance function $C(r,s,\theta
_{0})=\lim_{n\rightarrow \infty }C_{n}(r,s,\theta _{0})$.  
The following proposition establishes the limiting distribution of our
test statistics.

\begin{prop}
\label{propnull}Suppose Assumptions 3.1-3.3 hold. Then under $H_{0}$
\begin{equation*}
\Gamma (\hat{V}_{2n}(r))\stackrel{d}{\rightarrow}\Gamma (\hat{V}_{2\infty }(r)),
\end{equation*}%
where
\begin{equation*}
\hat{V}_{2\infty }({r})=V_{2\infty }(r)-h(r)^{\prime }\psi _{\infty }.
\end{equation*}
\end{prop}

The next
proposition provides results on asymptotic distribution  under the local alternatives.

\begin{as}
{The conditional cdfs $H_{t}(y|\Omega _{t})$ are continuous and
strictly increasing in $y$.}
\end{as}

\begin{prop}
\label{propalt} Suppose Assumptions 3.1-3.4 hold. Then under $H_{1n}$
\begin{equation*}
\Gamma (\hat{V}_{2n}(r))\stackrel{d}{\rightarrow }\Gamma (\hat{V}_{2\infty
}(r)+\delta k(r)-\delta h(r)^{\prime }\mu),
\end{equation*}%
where%
\begin{eqnarray*}
k(r) &=&{\plim_{n\to\infty}}\frac{1}{n}\sum_{t=2}^{n}\left\{ \left[
H_{t-1}(F_{t-1}^{-1}(r_{2}|\Omega _{t-1},\theta _{0})|\Omega _{t-1})-r_{2}%
\right] r_{1}\right. \\
&&\nsp\nsp\nsp\left. +\left[ H_{t}(F_{t}^{-1}(r_{1}|\Omega _{t},\theta _{0})|\Omega
_{t})-r_{1}\right]
I\left( U_{t-1}\leq H_{t-1}(F_{t-1}^{-1}(r_{2}|\Omega _{t-1},\theta _{0})|\Omega _{t-1})\right) \right\}. 
\end{eqnarray*}%
\end{prop}

Under $G_{nt}$,  $U_{t}\nsp=\nsp F_{t}(Y_{t}|\nsp\Omega _{t},\nsp\theta _{0})$ are not
 iid anymore; instead  $U_{t}^{\dag}\nsp=\nsp G_{nt}(Y_{t}|\nsp\Omega _{t},\nsp\theta _{0})$
are uniform iid. Due to this fact, we have a drift $k(r)$ in the asymptotic
distribution. Let us now examine  this drift more closer. If we fix $r_{2}=1$,
we get the drift of the process $\hat{V}_{1n}$ equal $k\left( \left(
r_{1},1\right) \right) =\limfunc{plim}\frac{1}{n}\sum_{t=2}^{n}\left[
H_{t}(F_{t}^{-1}(r_{1}|\Omega _{t},\theta _{0})|\Omega _{t})-r_{1}\right]$.
This drift might be zero even if $%
H_{t}$ and $F_{t}$ are different. If they differ only by the conditioning
set, averaging may smooth away this difference. The extreme is in the case
of elliptical distribution, where after integrating out one variable from
the conditioning set we are still in the same family, so the drift is zero.
As an example, consider testing the AR(1) model against AR(2),
both with standard normal innovations. This is equivalent to testing $H_{0}:$ 
$Y_{t}|\Omega _{t}$ $\sim $ $N(\alpha Y_{t-1},\sigma ^{2})$, for some $\alpha
$ and $\sigma $, which may be consistently estimated under $H_{0}$ by $\hat{%
\alpha}=\sum_{t=1}^{n}Y_{t}Y_{t-1}/\sum_{t=1}^{n}Y_{t}^{2}$ and $\hat{\sigma}%
^{2}=\frac{1}{n}\sum_{t=1}^{n}(Y_{t}-\hat{\alpha}Y_{t-1})^{2}$ (to consider
simple hypothesis we may fix parameters in $H_{0}$ to $%
\hat{\alpha}$ and $\hat{\sigma}$). In other words,
\begin{equation*}
(Y_{t}-\alpha Y_{t-1})/\sigma =:\varepsilon _{t}\sim \text{ iid with
cdf }F_{\varepsilon }(\varepsilon ).
\end{equation*}%
In this example $U_{t}=F_{\varepsilon }((Y_{t}-\alpha Y_{t-1})/\sigma)$.
Assume now that the true data generating process (DGP) is $Y_{t}|\Omega _{t}$ $\sim $ $%
N(\alpha _{1}Y_{t-1}+\alpha _{2}Y_{t-2},\sigma_2 ^{2})$. 
As before, denote the true distribution
$H_{t}(\cdot )$ and the null distribution $F_{t}(\cdot )$
(with a small abuse of notation since both depend on $%
\Omega _{t}$). We first show that $U_{t}$ are uniform (but not independent),
hence the  unconditional expectation of $V_{1n}$ is zero. Indeed, for some $\alpha_1,\alpha_2,\sigma_2^{2}$,%
\begin{equation*}
P(U_{t}\leq r)=
E\left[E\left[ I(Y_{t}\leq
F_{t}^{-1}(r))|\Omega _{t}\right]\right] =E\left[E\left[ H_{t}\left(F_{t}^{-1}(r)\right)|Y_{t}\right]\right]
=r,
\end{equation*}%
where in the last equality we use the particular form of $F_t(\cdot )$ and $%
H_{t}(\cdot )$ and the property of Gaussian distribution that $E\left[
H_{t}(\cdot )|Y_{t}\right] =F_{t}(\cdot )$ for chosen parameters. Therefore $E\left[V_{1n}(r)\right]= 0$. %

The drift $k(r)$
for $\hat{V}_{2n}$ can be written as $k\left( \left( r_{1},1\right) \right)
r_{2}+k\left( \left( 1,r_{2}\right) \right) r_{1}$ $+k_{2}\left( r\right) $,
where $k_{2}(r)$ is%
\begin{equation*}
\plim_{n\to\infty}\frac{1}{n}\nsp\sum_{t=2}^{n}\nsp\left(
H_{t}(F_{t}^{-1}(r_{1}|\Omega _{t},\nsp\theta _{0})|\Omega _{t})\nsp-\nsp r_{1}\right)
\nsp\left( I\left( U_{t-1}\nsp\leq\nsp H_{t-1}(F_{t-1}^{-1}(r_{2}|\Omega _{t-1}\nsp,\nsp\theta _{0})|\Omega _{t-1})\right) 
\nsp-\nsp r_{2}\right) .
\end{equation*}%
The term $k_{2}(r)$ gives additional power, by preventing the averaging out of
dynamic misspecification. In the aforementioned extreme case of dynamic
misspecification of elliptical distribution, this term alone delivers the
power reported in our Monte Carlo simulations. In terms of $U_{t},$ $%
k_{2}(r) $ controls the dynamics of $U_{t},$ while the other terms control
uniformity.

Under the alternative, we may have also left (\ref{eqestlinear}) not centered, then  $\mu\ne 0$. This term does not appear in methods which use projections, as in Bai (2003). 

Now, we discuss the consistency of the test against the fixed alternative.

\begin{as}
The following limit in probability exists for $r\in [0,1]^2$:
\begin{equation*}
P_2\left(r\right)=\plim_{n\to\infty} \frac{1}{n}\sum^n_{t=1}I\left(U_t\le r_1\right)I\left(U_{t-1}\le r_2\right),\quad\forall\theta_0\in\Theta.
\end{equation*}
\end{as}

Assumption 3.5 holds 
under stationarity and ergodicity of $U_t$, in which case $P_2\left(r\right)=P\left(U_t\le r_1,U_{t-1}\le r_2\right)$ for all $t$ whenever we are under the null or not. Moreover, under the null $P\left(U_t\le r_1,U_{t-1}\le r_2\right)= r_1 r_2.$ 
Let $\bar P_2(r)=P_2\left(r\right)- r_1 r_2$. Define the fixed alternative\footnote{We are not aware of any published paper with a consistent test for $H_0$. According to the working paper of Bierens and Wang (2014), in case the  assumption of strict stationarity is maintained the only available consistent test for $H_0$  is their test, which is based on comparing model-implied and model-free estimates of conditional characteristic functions and uses the approach of Bierens (1984) to deal with conditioning sets. Our procedure, in contrast,  is not restricted to strictly stationary data and is based on the probability integral transform, which delivers standard uniform and independent random variables both for stationary and nonstationary distributions and is already widely used in practice, as we discussed in the Introduction.}:

\bigskip

$H_{1}$: The conditional distribution function of $Y_{t}$ is equal to $H_{t}(y|\Omega _{t}),$ 
which is different from the assumed distribution $F_{t}(Y_t|\Omega _{t},\theta )$ in the following sense:  
\begin{equation*}
\forall\theta_0\in \Theta, \ \exists r\in [0,1]^2, \quad\text{ s.t. } \bar P_2(r)\ne 0.
\end{equation*}

\bigskip 
Under $H_1$, $V_{2n}(r)={n}\bar P_2(r)/{\sqrt{n-1}}$ diverges at least for some $r$, and therefore the KS test $\sup_r|V_{2n}(r)|$ is consistent. For the CvM test, we need  $\bar P_2(r)\ne 0$ at the set of $r$ with a positive measure. It might be the case that the null is violated, say $P\left(U_t\le r_1,U_{t-1}\le r_2\right)\ne r_1 r_2$, but $P\left(U_t\le r_1\right)= r_1$ (see the example above). In this case we can distinguish the alternative with $V_{2n}$ but not with $V_{1n}$. It might be also the case that the null is violated, say $P\left(U_t\le r_1,U_{t-2}\le r_2\right)\ne r_1 r_2$, but $P\left(U_t\le r_1,U_{t-1}\le r_2\right)= r_1 r_2.$ Therefore a test based solely on $V_{2n}$ is not consistent against the whole complement of the $H_0$ (unless time series are restricted to have first-order dependence only). Tests based on a combination of $V_{2n,j}$ for different $j$ are consistent against a broader set of alternatives: they check marginals and pairwise structure of $U_t$, but not serial structure. In theory we need to aggregate all $V_{2n,j}$ and $V_{pn}$, but this approach might not work well in finite samples. We now state the formal result.

\begin{as}
{$\hat{\theta}\stackrel{p}{\rightarrow }\theta _{1}$ for some $%
\theta _{1}\in \Theta .$}
\end{as}

Under an alternative, the estimator must converge in probability.
Under the null together with Assumption 3.3, this would imply $%
\theta _{1}=\theta _{0}$. Otherwise this is not necessarily true,
and $\theta _{1}$ is often called a ``pseudo-true" value.

\begin{prop}
\label{propFA}
Suppose Assumptions 3.1, 3.2, 3.4-3.6 hold. Then under $H_1$ for all sequences of rv's $c_n=O_p(1)$ we have
\begin{equation*}\label{eq:GVFA}
\lim_{n\to\infty}P\left(\Gamma (\hat{V}_{2n}({r}))>c_n\right)=1.
\end{equation*}
\end{prop}

We show that bootstrap critical values are bounded both under the null and under alternative. Therefore Proposition  \ref{propFA} suffices for the consistency of the bootstrap assisted test. We prove that $\Gamma (\hat{V}_{2n}^{\ast })$ has the same limiting
distribution as $\Gamma (\hat{V}_{2n})$. The proof is similar to Andrews (1997); we need to establish an analog of his (4.3) for introduced tests and in the dynamic setup. 
We say that the sample is distributed ``\textit{under }$\{\theta _{n}:n\geq
1\} $'' when there is a triangular array of rv's $\{Y_{nt}:n\geq 1,t\leq n\}$
with $(n,t)$ element generated by $F_{t}(\cdot |\Omega _{nt},\theta _{n})$,
where $\Omega _{nt}=(Y_{nt-1},Y_{nt-2},...)$.

\begin{as} 
{
\label{asslb} For all nonrandom sequences $\{\theta _{n}:n\geq 1\}$ for
which $\theta _{n}\rightarrow \theta _{0}$, we have
\begin{equation}\label{eq:linexptr}
\sqrt{n}(\hat{\theta}-\theta _{n})=\frac{1}{\sqrt{n}}\sum_{t=1}^{n}\psi\left(Y_{nt},\Omega_{nt},\theta_n\right)+o_{p}(1),
\end{equation}%
under $\{\theta _{n}:n\geq 1\}$, where 
\begin{equation*}
\frac{1}{n}\sum_{t=1}^{n}\psi\left(Y_{nt},\Omega_{nt},\theta_n\right)\psi\left(Y_{nt},\Omega_{nt},\theta_n\right)^{\prime
}\stackrel{p}{\rightarrow }\Psi .
\end{equation*}}
\end{as}

Note that the function $\psi$ is the same as in
Assumption 3.3. We require that estimators of values of $\theta$ close to $\theta _{0}$
have \emph{the same} linear representation as the estimator of $%
\theta _{0}$ itself.  Assumption 3.7 is not much more restrictive than Assumption 3.3, since most proofs of linear expansion of parametric estimators can be accommodated for the triangle linear expansion (\ref{eq:linexptr}). The next proposition states that the
asymptotic distribution of the test statistics with the
bootstrapped data, or denoted shortly under $\{\theta _{n}:n\geq
1\} $,  coincides with the prior result obtained under the null.

\begin{prop}
\label{propboot}Suppose Assumptions 3.1, 3.2 and 3.7 hold. Then under $H_{1n}$, for any nonrandom sequence $%
\{\theta _{n}:n\geq 1\}$ for which $\theta _{n}\rightarrow \theta _{0}$,
under $\{\theta _{n}:n\geq 1\}$%
\begin{equation}\label{eq:GVB}
\Gamma (\hat{V}_{2n}({r}))\stackrel{d}{\rightarrow }\Gamma (\hat{V}_{2\infty
}({r})).
\end{equation}
\end{prop}

Let $c_{\alpha n}\left(\theta_n\right)$ denote the level $\alpha$ critical value of $\Gamma (\hat{V}_{2n}({r}))$ generated with some $\theta_n$. Let $c_{\alpha }\left(\theta_0\right)$ denote the level $\alpha$ critical value of $\Gamma (\hat{V}_{2\infty}({r}))$. By Proposition~\ref{propboot} and absolute continuity of limiting distribution  $\Gamma (\hat{V}_{2\infty}({r}))$, $c_{\alpha n}\left(\theta_n\right)\to c_{\alpha }\left(\theta_0\right)$ with probability $1$.
Then, if $\hat\theta\stackrel{p}{\rightarrow }\theta_1$, $c_{\alpha n}\left(\hat\theta\right)\stackrel{p}{\rightarrow } c_{\alpha }\left(\theta_1\right)$. Randomness here comes both from the sample and the bootstrap simulations. Then the asymptotic significance level of the $\Gamma (\hat{V}_{2n}({r}))$ with critical value $c_{\alpha n}\left(\hat\theta\right)$ is $\alpha$. $c_{\alpha n}\left(\hat\theta\right)$ in turn, is approximated by $c_{\alpha n B}\left(\hat\theta\right)$ when $B\to\infty,$ where the latter is the  critical value after $B$ repetitions. Note, that we need only $\hat\theta\stackrel{p}{\rightarrow }\theta_1$, so the analysis holds under the null and under the alternative (in the latter case with Assumption 3.6 and Assumptions 3.2 and 3.7 holding for any $\theta_0\in\Theta$), thus justifying the bootstrap approximation.
For more details, see the discussion on pages 1107-1109 of Andrews (1997).

\section{FINITE SAMPLE PERFORMANCE AND EMPIRICAL APPLICATION}
\setcounter{theorem}{0}\setcounter{equation}{0}
In this section, we report the results of a Monte Carlo study to investigate the
finite sample performance of the proposed  tests.
The number of Monte Carlo repetitions is set to $1000$. To calculate critical values we use the fast bootstrap method of Giacomini, Politis and White (2013). To save space we present here in detail a case of conditional mean misspecification in GARCH models, while other simulation results can be obtained upon request. We also show how to apply our technique to  models of stock exchange indexes. The models are estimated by the method of Maximum Likelihood (ML). For GARCH processes we consider a stationary solution, therefore the ULLN required for Assumption~3.2 holds, and the ML estimator is consistent and satisfies Assumption~3.3; see Fan and Yao (2003) for conditions on GARCH and examples of other nonlinear time series models which deliver stationary and ergodic/mixing solutions.

While our methodology applies in very general contexts and there are potentially more powerful tests specifically designed for GARCH models\footnote{For known specific alternatives, Jarque-Bera, Box-Pierce, and Engle's test for ARCH effects and other parametric tests are more powerful than nonparametric tests in general and those proposed here in particular. However, they may have no power against other alternatives. Since the goal of the paper is to test the null against a general alternative, parametric tests are not included in the  comparison.}, we stick to the GARCH-type null models for  two reasons. First,  we make evident  the ability of our test to detect misspecification both in dynamics and marginal distributions, which are easy to introduce to these models and which are usually tested separately. Second, the combination of the probability integral transform and residual empirical process is employed in  Bai (2003), who applies his tests to GARCH models. We discussed above the limitations of (univariate) residual empirical process tests and our goal is to verify that our approach indeed outperforms them.  Thus, we compare  the performance of our test with the ones based on a univariate process for the same null hypothesis and the same data. In order to make tests comparable, we avoid the martingale transform and use bootstrapped versions of all tests. Our methodology applies without modification to all numerous extensions of the basic GARCH model, including models with asymmetry, leverage, higher moments and nonlinear moment dynamics. 

\subsection{Conditional mean specification in GARCH models}
\begin{figure}[!h]
\centering 
\subfigure[Tests $D_{1,100}^{CvM}$ and $ADJ_{j,100}^{0}$] {
\includegraphics[width=0.47\textwidth]{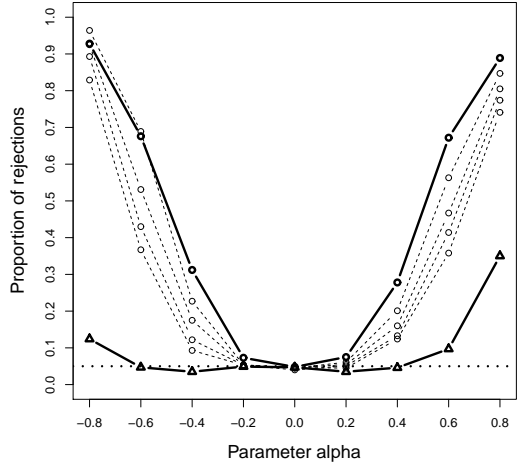}
}
\subfigure[Tests $D_{1,100}^{CvM}$ and $ADJ_{j,100}$] {
\includegraphics[width=0.47\textwidth]{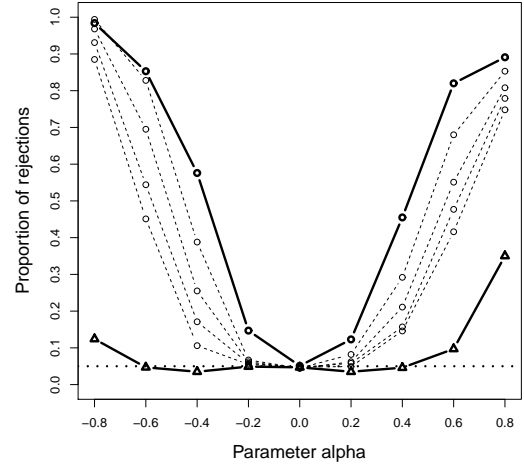}
}
\subfigure[Tests $D_{1,300}^{CvM}$ and $ADJ_{j,300}^{0}$] {
\includegraphics[width=0.47\textwidth]{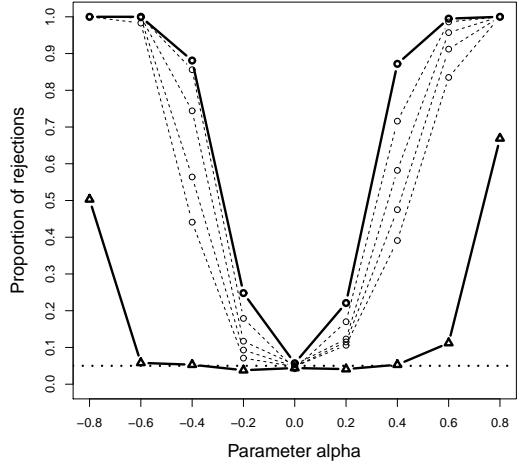}
}
\subfigure[Tests $D_{1,300}^{CvM}$ and $ADJ_{j,300}$] {
\includegraphics[width=0.47\textwidth]{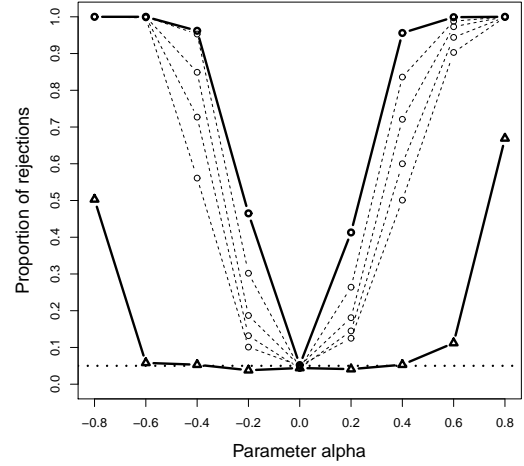}
} 
\caption{Proportion of rejections of tests for the GARCH(1,1) model and DGP AR(1)-GARCH(1,1), given in (\ref{eq:ar10g}). In both models, innovations are independent Gaussian.
On Panels (a) and (c), tests based on a one-parameter process $D_{1n}^{CvM}$ (thick line with triangles markers) and two-parameter process with 1 lag $ADJ_{1n}^{0}$ (thick line with circles markers) are considered, as well as two-parameter processes with $j\in\{2,3,4,5\}$ lags $ADJ_{jn}^{0}$ (dashed lines with circles markers). On Panels (b) and (d) $ADJ_{jn}^{0}$ are changed to $ADJ_{jn}$. Rejections at $5\%$ are plotted with a dotted line.
Sample sizes are $n=100$ on Panels (a) and~(b); $n=300$ on Panels (c) and (d).}
\label{fig:ar10garch}
\end{figure}

In this experiment, we examine the GARCH(1,1) model against the AR(1)-GARCH(1,1) data generating process (DGP).
This example is motivated by the findings in our application for the tests on real data (see below). The null models are GARCH(1,1), and DGP is
\begin{equation}\label{eq:ar10g}
Y_t=\alpha_1 Y_{t-1}+h_t\varepsilon_{t}
\end{equation}
with
\begin{equation*}
h^2_t=0.1+0.1 \left(Y_{t-1}-\alpha_1 Y_{t-2}\right)^2+0.8h^2_{t-1}.
\end{equation*}
Parameter $\alpha_1$ takes values $-0.8, -0.6,...,0,...,0.6,0.8$; innovations $\varepsilon_t$ in both models are independent Gaussian, sample sizes $n=100$ and $n=300$.
In Figure~\ref{fig:ar10garch}, 
\begin{figure}[!h]
\centering 
\subfigure[Tests $D_{1,100}^{CvM}$ and $ADJ_{j,100}^{0}$] {
\includegraphics[width=0.47\textwidth]{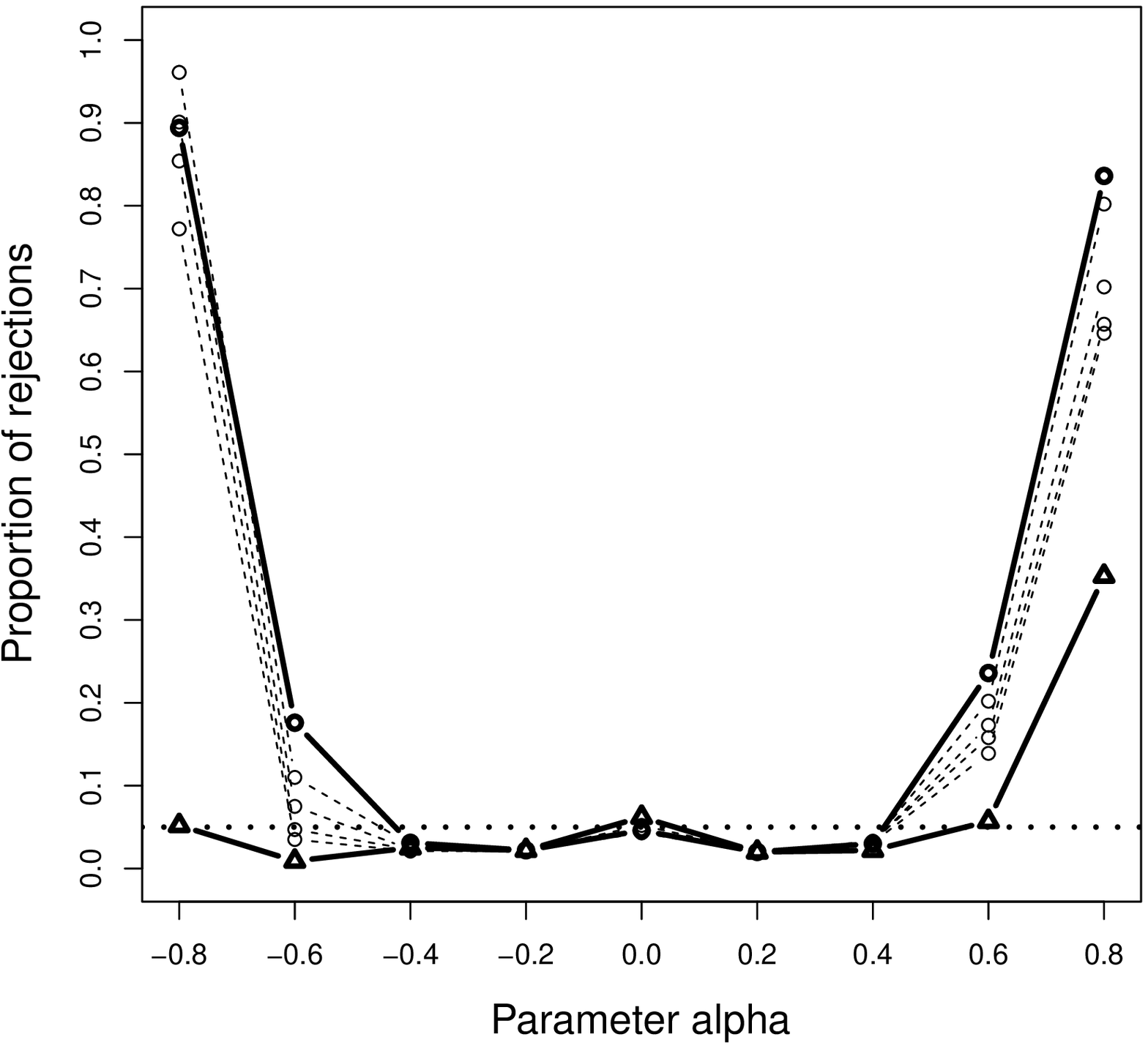}
}
\subfigure[Tests $D_{1,100}^{CvM}$ and $ADJ_{j,100}$] {
\includegraphics[width=0.47\textwidth]{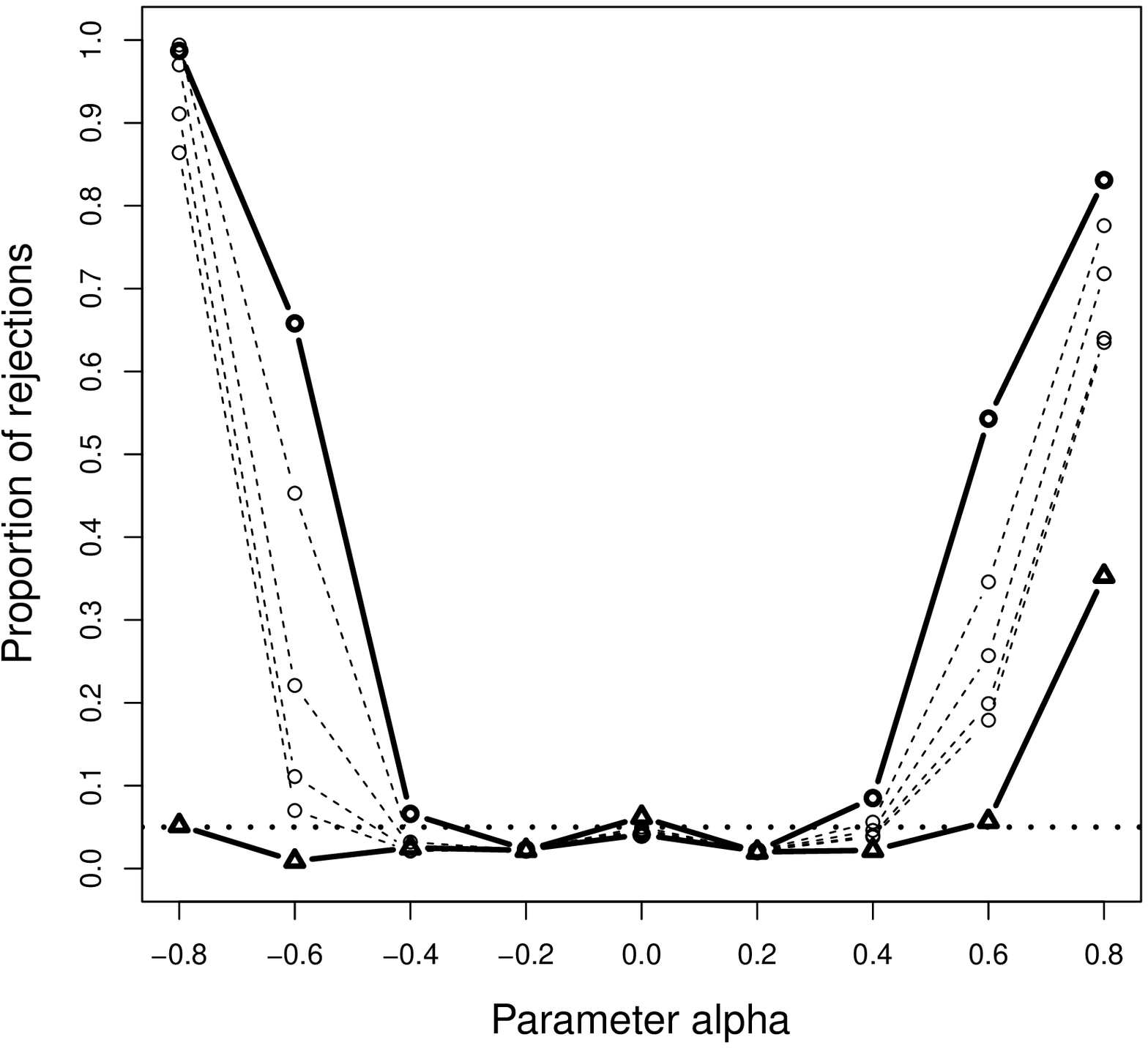}
}
\subfigure[Tests $D_{1,300}^{CvM}$ and $ADJ_{j,300}^{0}$] {
\includegraphics[width=0.47\textwidth]{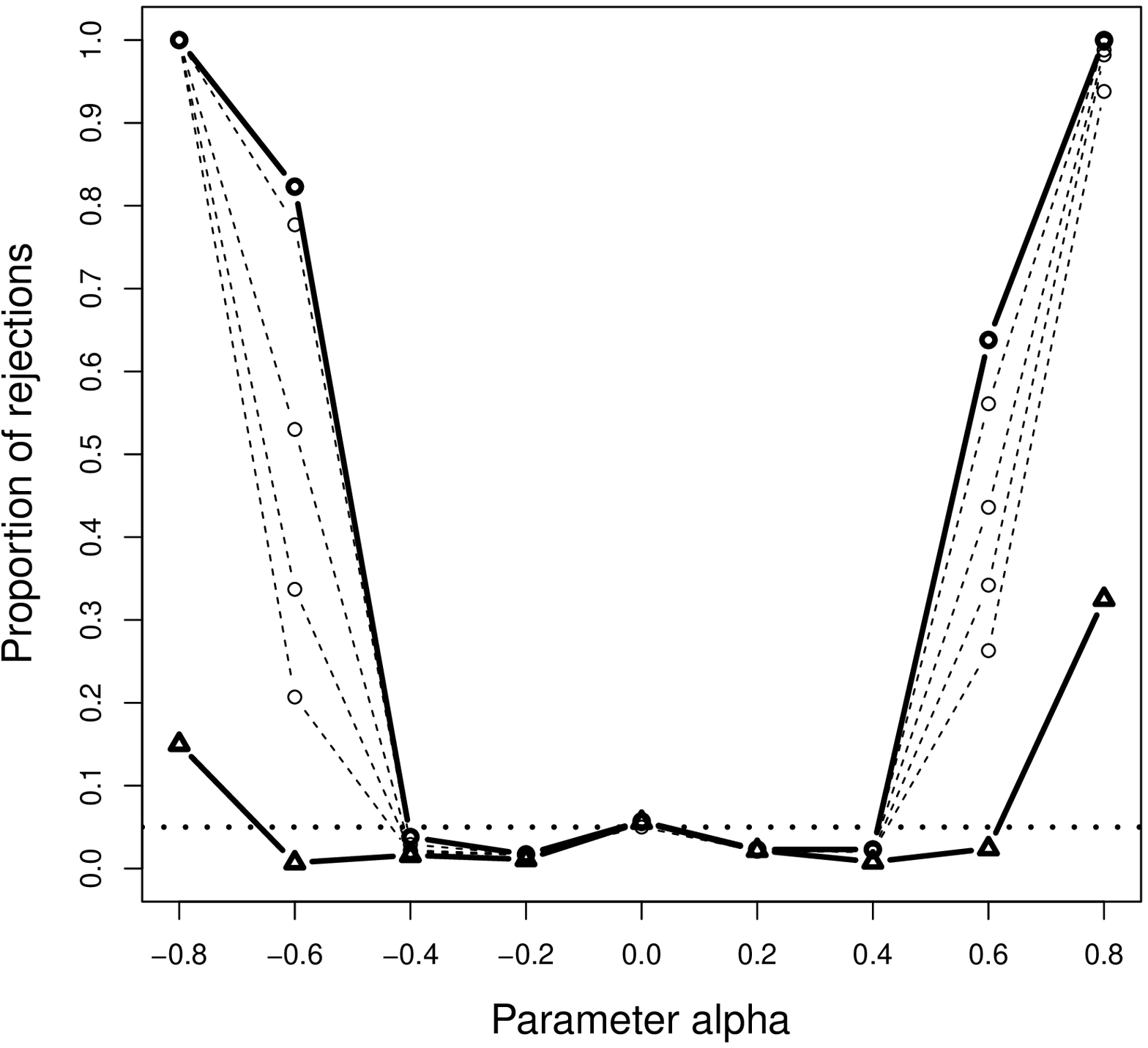}
}
\subfigure[Tests $D_{1,300}^{CvM}$ and $ADJ_{j,300}$] {
\includegraphics[width=0.47\textwidth]{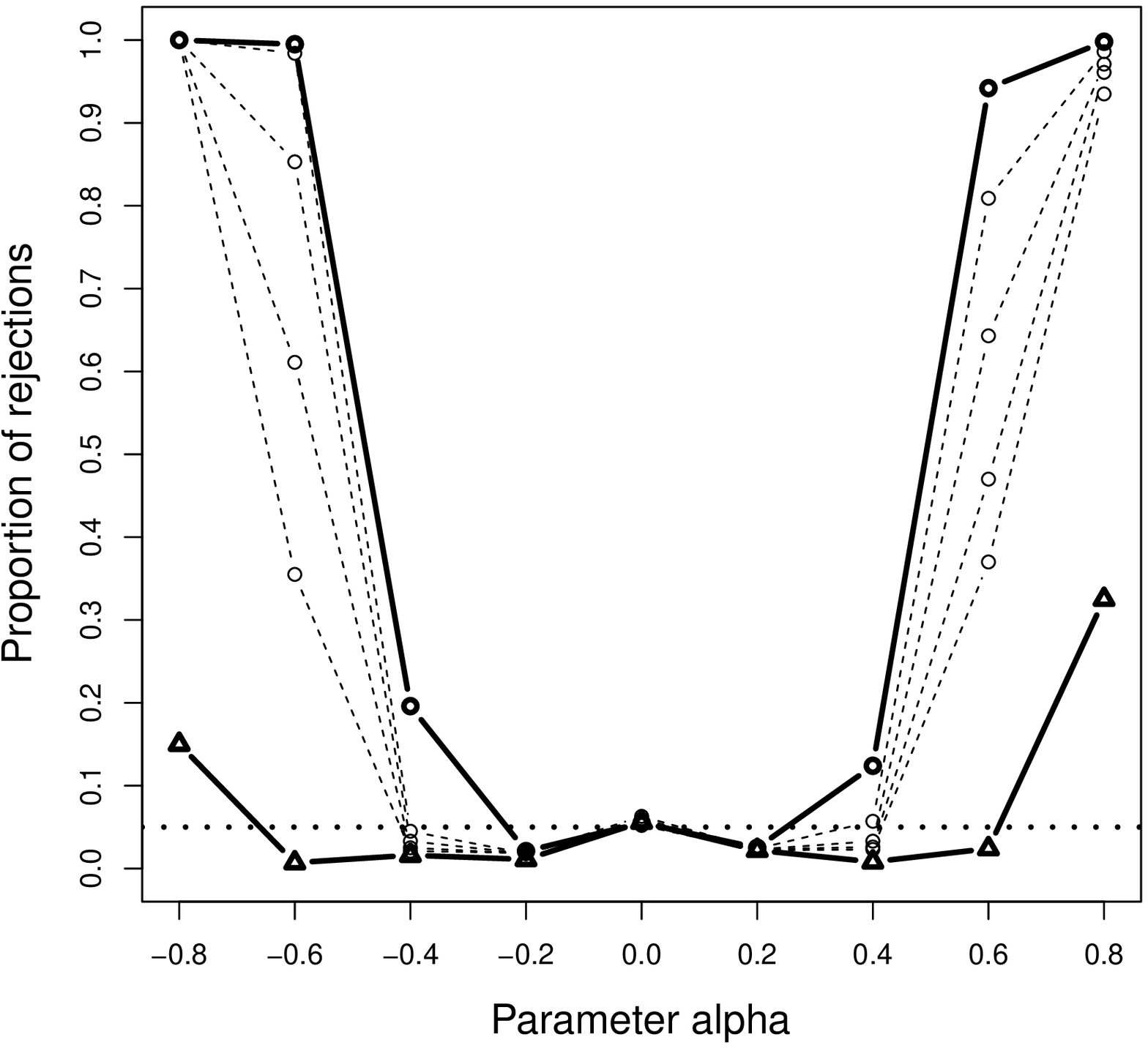}
} 
\caption{Proportion of rejections of tests for the GARCH(1,1) model and DGP AR(1)-GARCH(1,1), given in (\ref{eq:ar10g}). In both models, innovations are student-t with 5 degrees of freedom.
On Panels (a) and (c), tests based on one-parameter process $D_{1n}^{CvM}$ (thick line with triangles markers) and two-parameter process with 1 lag $ADJ_{1n}^{0}$ (thick line with circles markers) are considered, as well as two-parameter processes with $j\in\{2,3,4,5\}$ lags $ADJ_{jn}^{0}$ (dashed lines with circles markers). On Panels (b) and (d) $ADJ_{jn}^{0}$ are changed to $ADJ_{jn}$. Rejections at $5\%$ are plotted with a dotted line.
Sample sizes are $n=100$ on Panels~(a) and (b); $n=300$ on Panels (c) and (d).}
\label{fig:ar10garcht5}
\end{figure}
the top left plot provides the proportion of rejections at the $5\%$ level of tests $D^{CvM}_{1n}$ and $ADJ^0_{1n}$  against different parameters $\alpha_1$, with sample size $n=100$.
When $\alpha_1=0,$
we see that the size of  both tests is very close to a nominal $5\%$. The test $D^{CvM}_{1n}$ does not have power on the interval $[-0.6, 0.4]$ and has very low power for other parameter values. The test $ADJ^0_{1n}$ has power against all alternatives.
We also show results for $ADJ^0_{jn}$ (dashed lines with circle markers) with $j=2,3,...,5$ on the same plot. They fill monotonically the space from $ADJ^0_{1n}$ to $D^{CvM}_{1n}$
 with $j=2$ closest to $ADJ^0_{1n}$. The performance of these tests is decreasing with $j$. Because the misspecification is in the first  lag, most of the power comes from $D^{CvM}_{2n,1}$, capturing dependence between $U_t$ and $U_{t-1}$, although in this case there is also dependence in further lags, for instance between $U_t$ and $U_{t-2}$, $U_{t-3}$, etc. The more tests we aggregate, the less weight is given to the first lag and less power we have but we can capture a wider set of alternatives. The same effect is observed with Ljung-Box tests. If we do not include powerless $D^{CvM}_{1n}$ into aggregation (Figure~\ref{fig:ar10garch}, top right), the performance is better. On the bottom plots, we repeat the experiment for $n=300$. Here we see much better power for our tests, in particular this is close to $1$ starting from $|\alpha_1|=0.4.$

 We also run this experiment with student-t with 5 degrees of freedom innovations $\varepsilon_t$ (Figure \ref{fig:ar10garcht5}). Our tests do not have power for $-0.4\le\alpha_1\le 0.2$ for sample size $n=100$ and have power for all parameter values,
 for sample size $n=300$, whereas $D^{CvM}_{1n}$ has no power for  $\alpha_1\le 0.6.$ for both sample sizes. Similar results were obtained for tests based on the Kolmogorov-Smirnov norm but are not reported here to save space.
To summarize,   the tests have a size close to nominal, a good power for $n=100$, which increases when we increase the sample size (to $n=300$). This result is robust for the distribution of the innovations, i.e. we have similar pictures for both normal and student innovations.

\subsection{Application to stock exchange index}

\begin{figure}[!h]
{\centering
\includegraphics[width=0.47\textwidth]{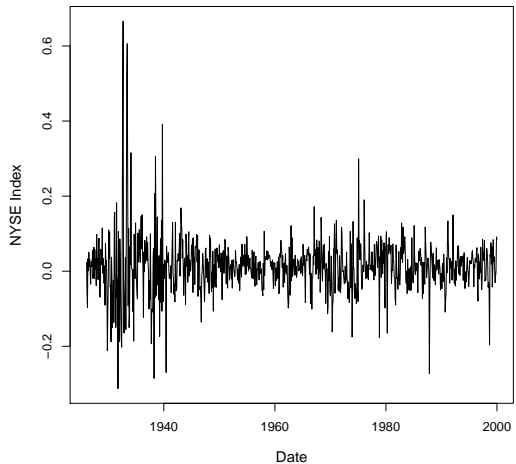}
\caption{Real data example: NYSE monthly equal-weighted returns 1926.1 - 1999.12.}\label{fig:data2}
}
\end{figure}

Consider the monthly NYSE equal-weighted returns for the data span from January 1926 to December 1999, see Figure \ref{fig:data2}.
Bai (2003) applied  transformed one-parameter test to this data, which rejected GARCH(1,1)-N at the 1\% significance level but could not reject GARCH(1,1)-$t_5$  (at 5\%).
In this example we illustrate how to use our ``generalized autocorrelations."

\begin{figure}[!h]
\centering 
\subfigure[CvM tests for GARCH(1,1)] {
\includegraphics[width=0.47\textwidth]{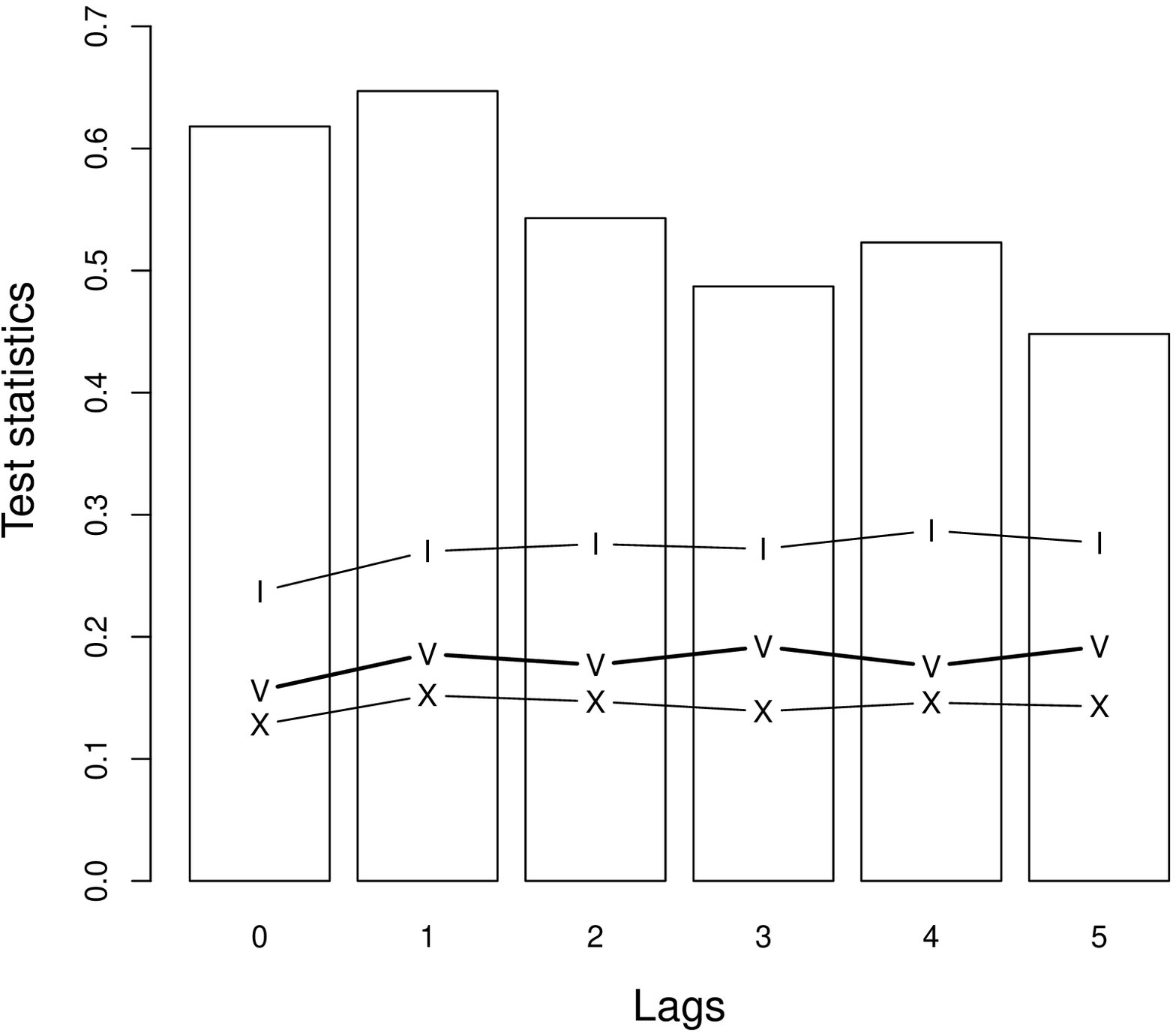}
}
\subfigure[KS tests for GARCH(1,1)] {
\includegraphics[width=0.47\textwidth]{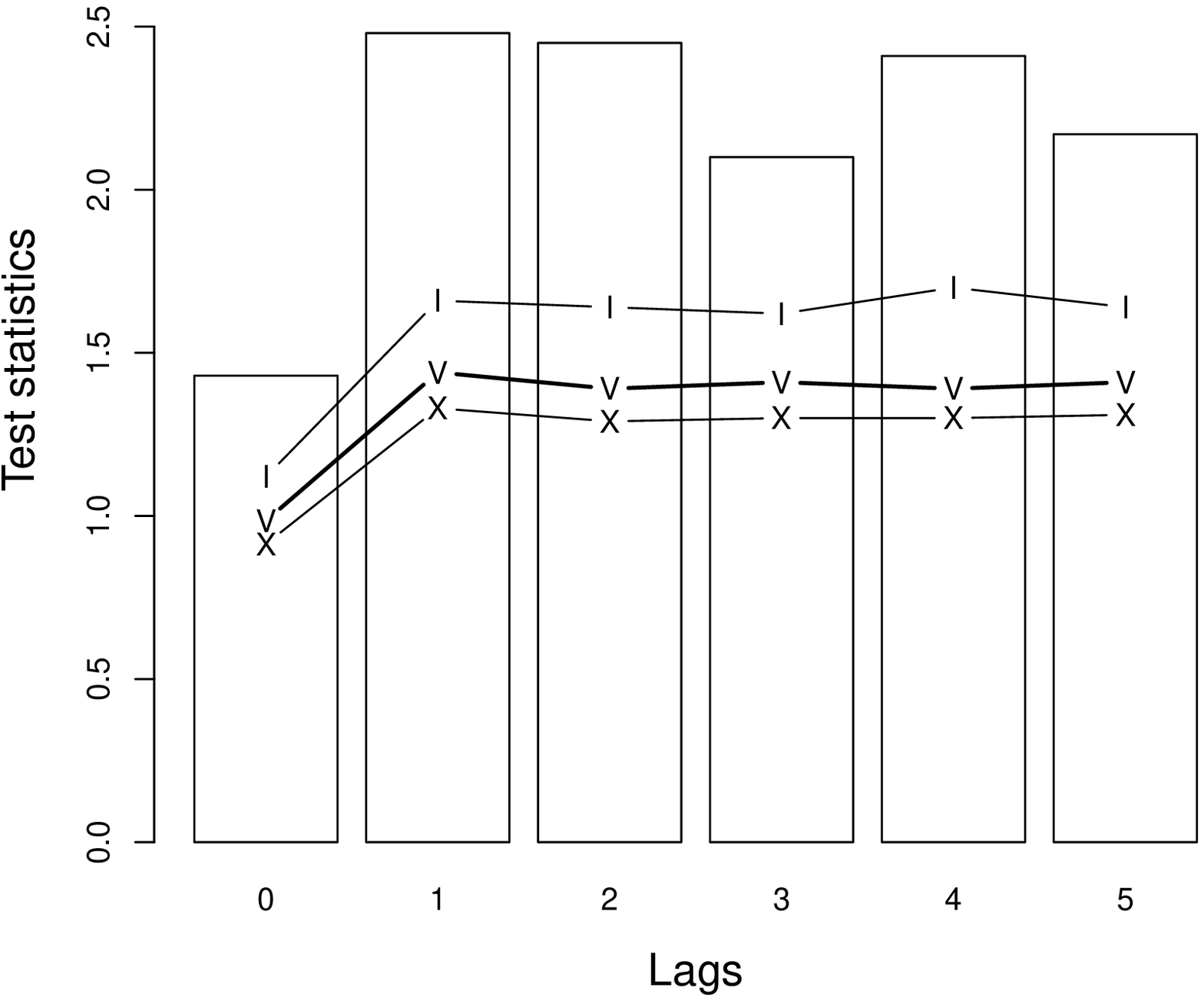}
}
\subfigure[CvM tests for GARCH(1,1)-$t_5$] {
\includegraphics[width=0.47\textwidth]{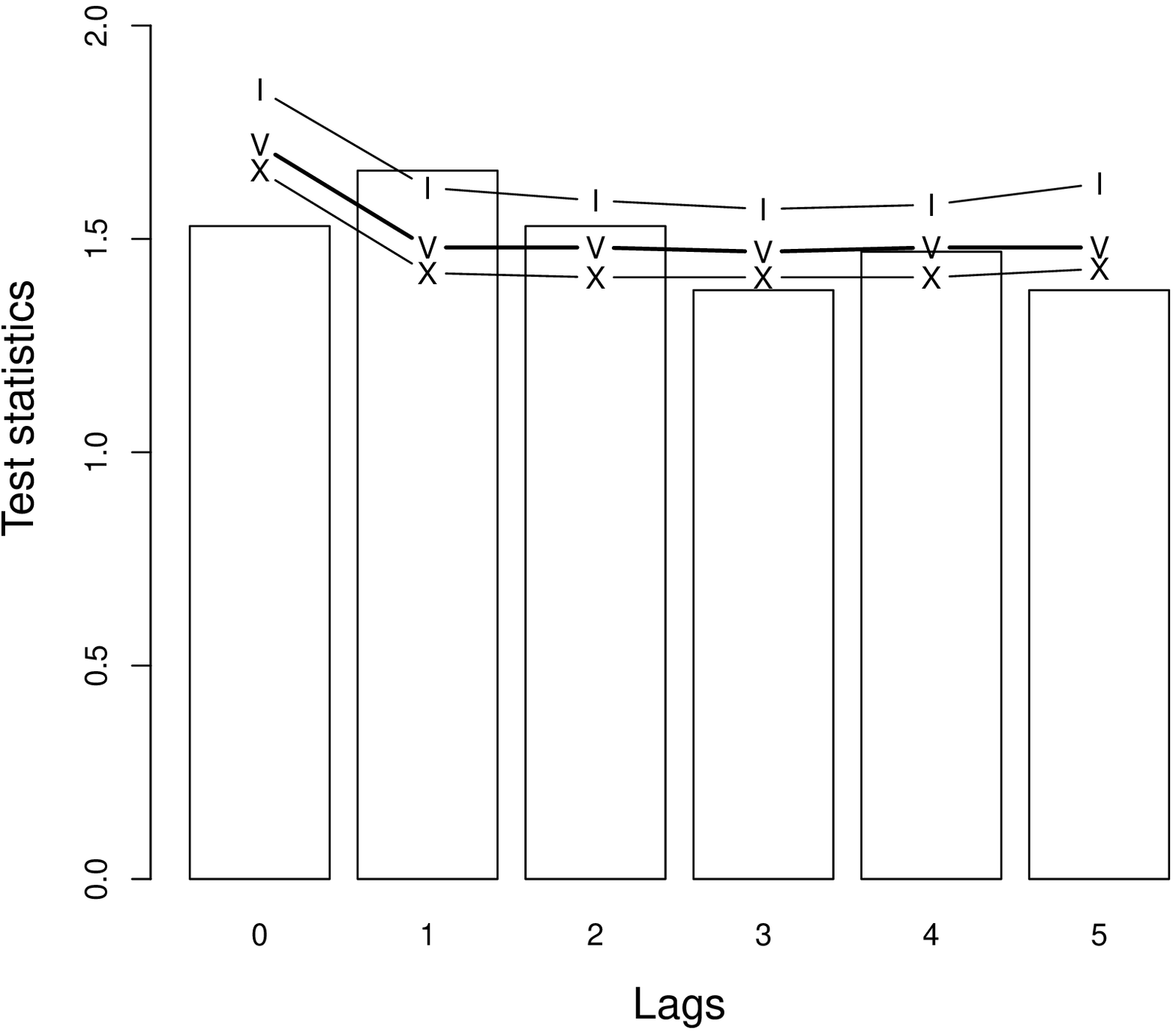}
}
\subfigure[KS tests for GARCH(1,1)-$t_5$] {
\includegraphics[width=0.47\textwidth]{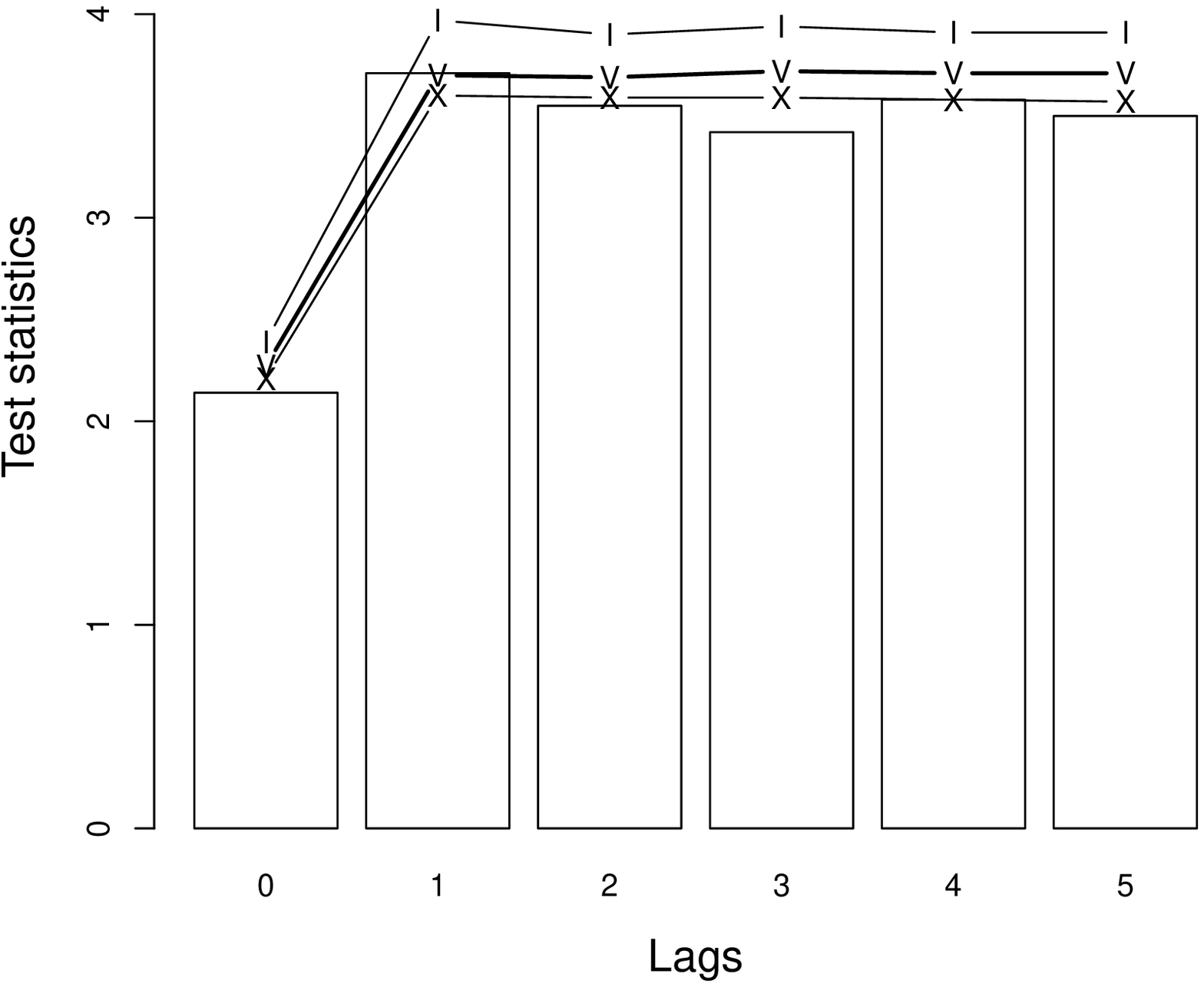}
} 
\caption{Real data example: testing models for NYSE monthly equal-weighted returns 1926.1 - 1999.12. Generalized autocorrelations (bars) and bootstrapped critical values (10\% - ``X", 5\% - ``V", 1\% - ``I") based on $D_{1n}$ (``lag $0$") and $D_{2n,j}$ (lags $j=1,2,3,4,5$) are plotted. Cramer-von Misses  and Kolmogorov-Smirnov tests for GARCH(1,1) model are shown on Panels (a) and (b) respectively. Cramer-von Misses  and Kolmogorov-Smirnov
 tests for GARCH(1,1)-$t_5$ model are shown on Panels (c) and (d) respectively.}
\label{fig:gauto-g11}
\end{figure}

We determine that $D_{1n}$ tests reject GARCH(1,1)-N at  the 1\% significance level but do not reject GARCH(1,1)-$t_5$ at the
10\% significance level.  At the same time bi-parameter tests reject both models (Table \ref{tab:baipval}, the first two lines). This is not surprising given our findings in Subsection 4.1.
If we check generalized autocorrelations (Figure \ref{fig:gauto-g11}), 
\begin{figure}[!h]
\centering 
\subfigure[CvM tests for AR(1)-GARCH(1,1)] {
\includegraphics[width=0.47\textwidth]{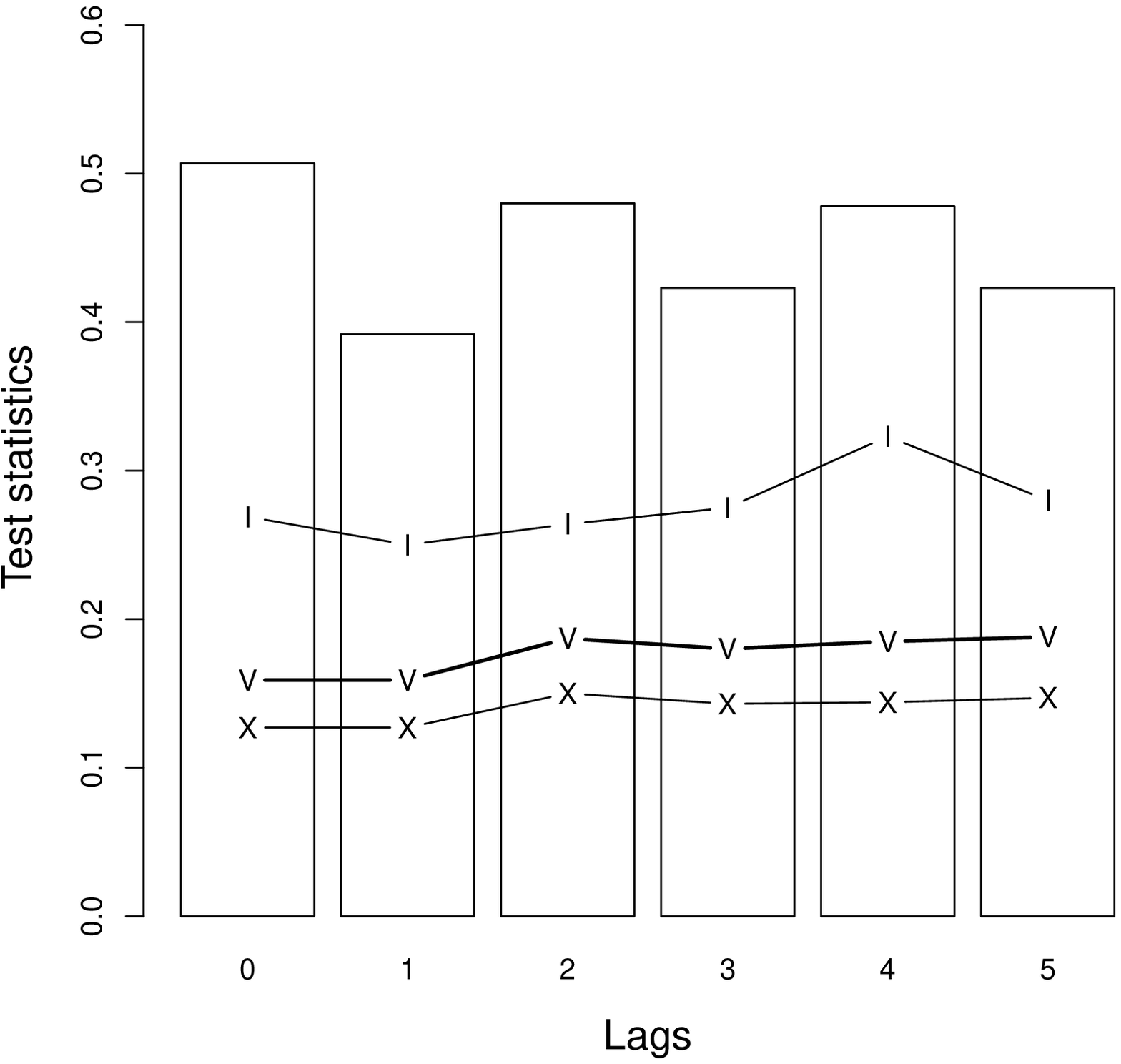}
}
\subfigure[KS tests for AR(1)-GARCH(1,1)] {
\includegraphics[width=0.47\textwidth]{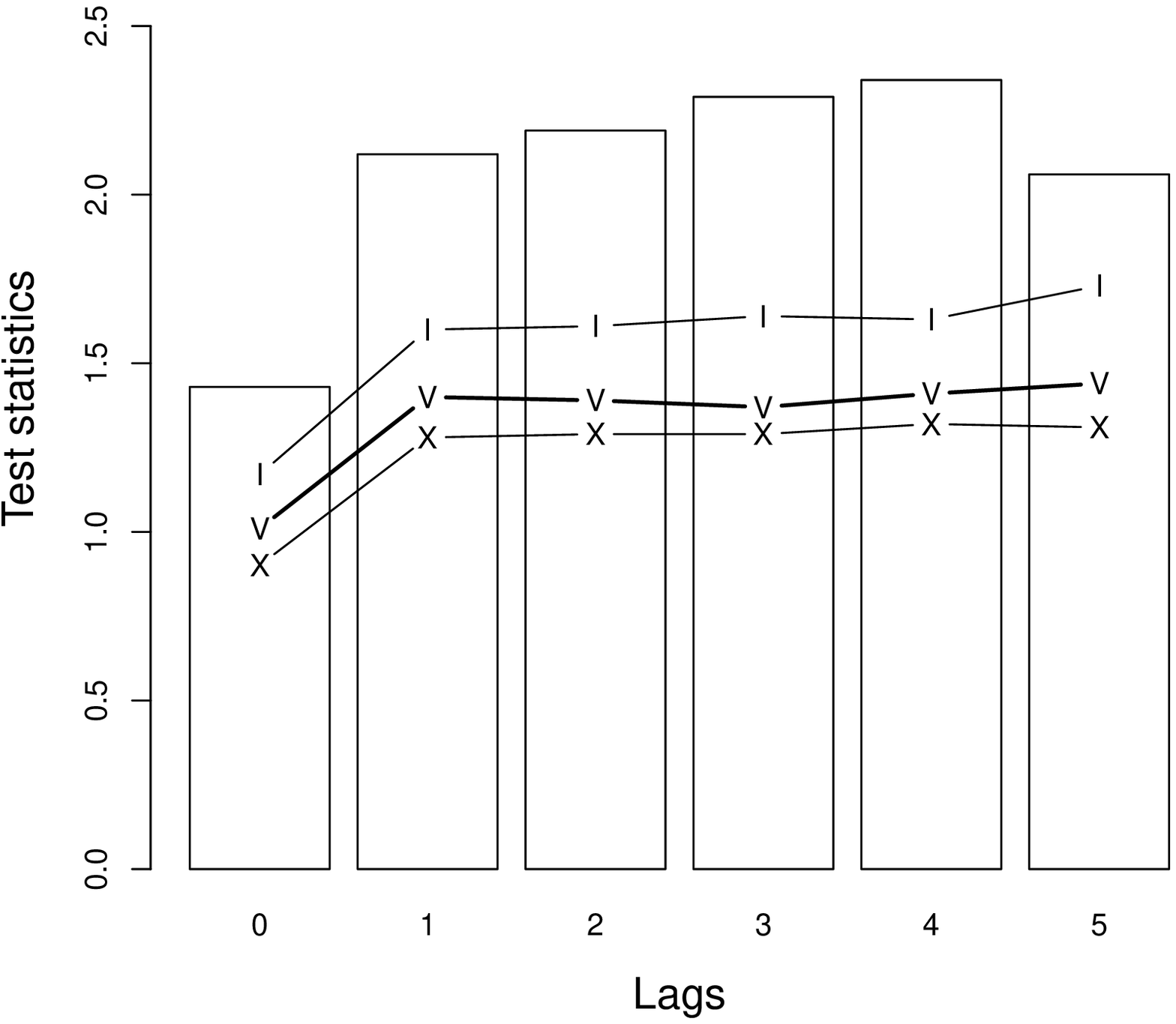}
}
\subfigure[CvM tests for AR(1)-GARCH(1,1)-$t_5$] {
\includegraphics[width=0.47\textwidth]{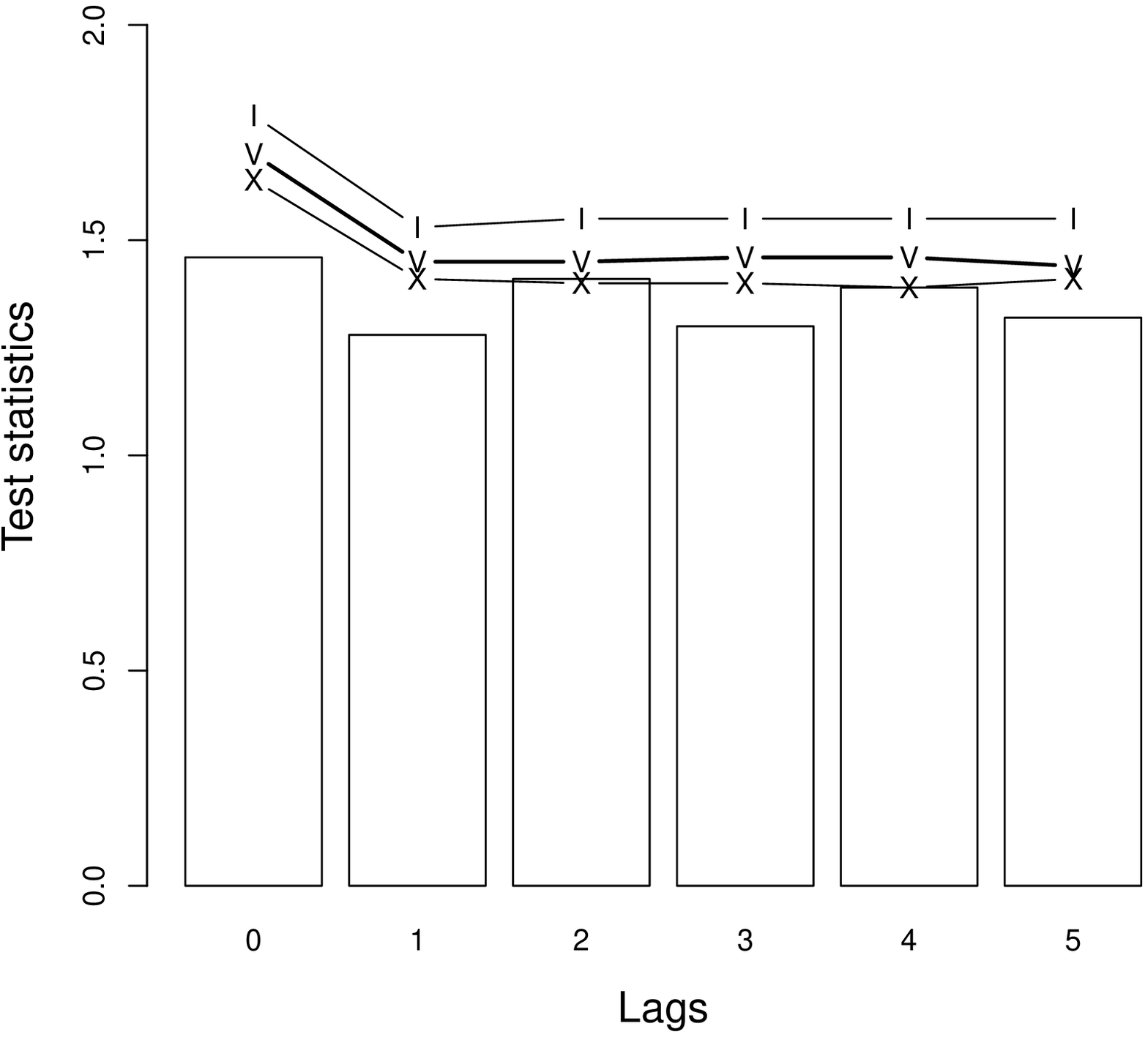}
}
\subfigure[KS tests for AR(1)-GARCH(1,1)-$t_5$] {
\includegraphics[width=0.47\textwidth]{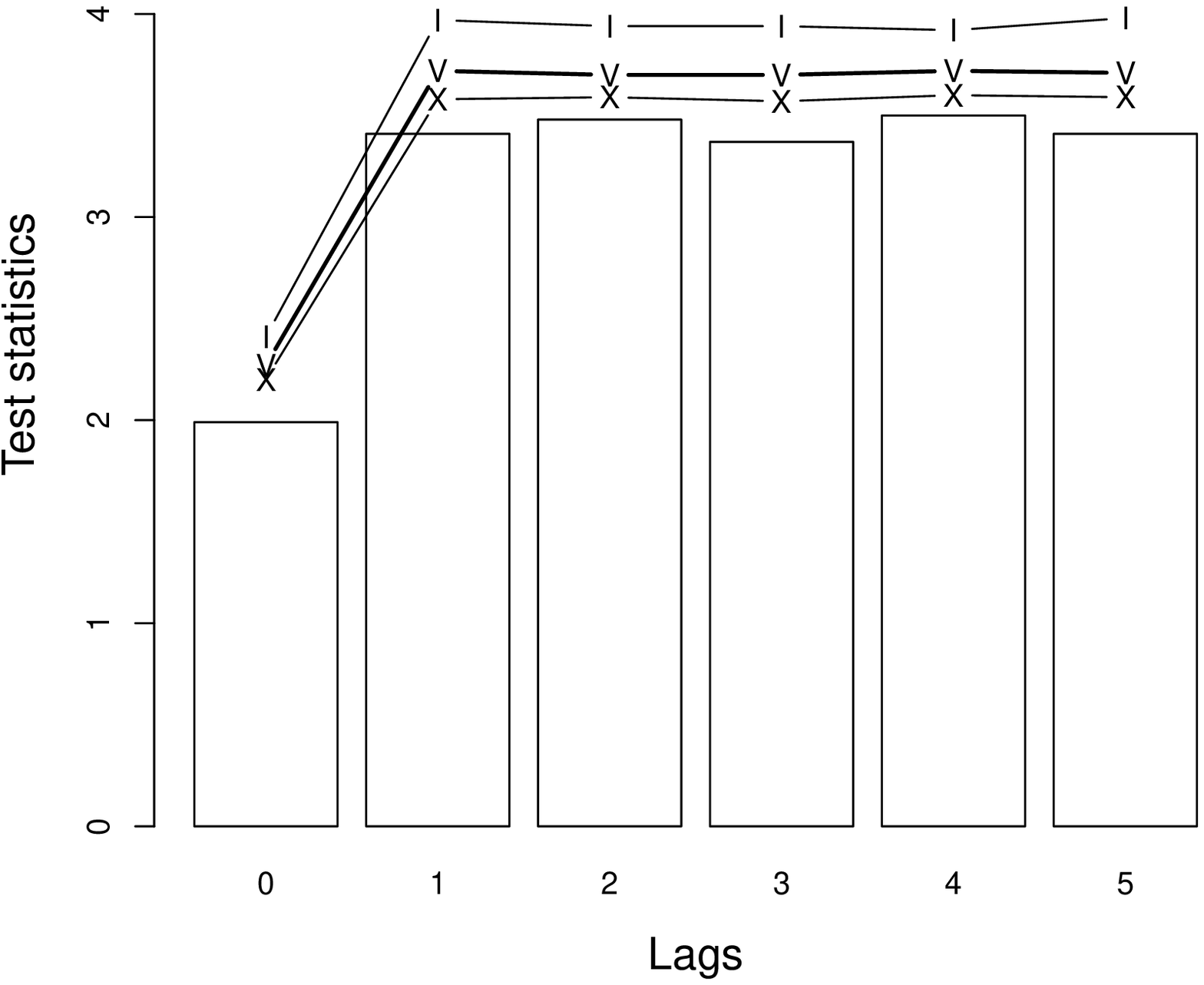}
} 
\caption{Real data example: testing models for NYSE monthly equal-weighted returns 1926.1 - 1999.12. Generalized autocorrelations (bars) and bootstrapped critical values (10\% - ``X", 5\% - ``V", 1\% - ``I") based on~$D_{1n}$~(``lag $0$") and $D_{2n,j}$ (lags $j=1,2,3,4,5$) are plotted. Cramer-von Misses  and Kolmogorov-Smirnov tests for AR(1)-GARCH(1,1) model are shown on Panels~(a) and (b) respectively. Cramer-von Misses  and Kolmogorov-Smirnov
 tests for AR(1)-GARCH(1,1)-$t_5$ model are shown on Panels~(c) and (d) respectively.}
\label{fig:gauto-ar1g11}
\end{figure}
we see that one-lag generalized autocorrelations are significant (at 1\% for CvM and at 5\% for KS) suggesting that not all dynamics is captured by the model.  One-parameter tests do not reveal it. To account for these dynamics we fit the AR(1)-GARCH(1,1) model. Generalized autocorrelations  are within the critical bounds for model with $t_5$ innovations (Figure \ref{fig:gauto-ar1g11}). All tests reject normal innovations while none of them  rejects $t_5$ at
1\% significance level.  See Table~\ref{tab:baiest} for ML estimates and their standard errors (in brackets) and Table~\ref{tab:baipval} for $p$-values of the test statistics for all considered models. We conclude that for the monthly NYSE equal-weighted returns for the data span from January 1926 to December 1999, the AR(1)-GARCH(1,1)-$t_5$ model cannot be rejected.

\begin{table}[tbp]
\caption{Real data example: estimates.}
{\centering
\begin{tabular}{@{}llccccc@{}}
\hline
 & & Mean const  & AR  & Variance const  & GARCH   & ARCH 
\tabularnewline
\hline
\hline
1 & GARCH(1,1) & $ 0.0130  $ & ---    & $ 0.0001 $ & $ 0.8433 $ & $ 0.1374 $    \tabularnewline
 &                       & $ (0.0018) $ &      & $ (0.00003) $ & $ (0.0209) $ & $ (0.0194) $   \tabularnewline
 \hline
2 & GARCH(1,1)-$t_5$ & $ 0.0142 $& ---    & $ 0.0002 $  & $ 0.8288 $  & $ 0.1459 $   \tabularnewline
& & $ (0.0016) $  &     & $ (0.00006) $  & $ (0.0344) $  & $ (0.0349) $    \tabularnewline
 \hline
3 & AR(1)-GARCH(1,1) & $ 0.0107 $  & $ 0.1939 $  & $ 0.0001 $  & $ 0.8468 $  & $ 0.1327 $    \tabularnewline
& & $ (0.0017) $  & $ (0.0367) $  & $ (0.00003) $  & $ (0.0219) $  & $ (0.0199) $    \tabularnewline
 \hline
4 & AR(1)-GARCH(1,1)-$t_5$  & $ 0.0117 $  & $ 0.1656 $  & $ 0.0001 $  & $ 0.8340 $  & $ 0.1417 $    \tabularnewline
& & $ (0.0016) $  & $ (0.0339) $  & $ (0.00001) $  & $ (0.0348) $  & $ (0.0351) $    \tabularnewline
 \hline
 \hline
\end{tabular}
}\label{tab:baiest}
\footnotesize
\renewcommand\baselineskip{11pt}
\textbf{Note:} Maximum likelihood estimates of GARCH-type models for NYSE monthly equal-weighted returns 1926.1 - 1999.12.
\end{table}

\begin{table}[tbp]
\caption{Real data example: $p$-values of specification tests.}
{\centering
\begin{tabular}{@{}llcccccc@{}}
\hline
 & $H_0$  & $ D_{1n}^{CvM} $  & $ ADJ_{1n} $  & $ ADJ_{5n} $   & $ D_{1n}^{KS} $  & $ MDJ_{1n} $  & $ MDJ_{5n} $
\tabularnewline
\hline
\hline
1 & GARCH(1,1) & $ 0.001^{***} $  & $ 0.001^{***} $  & $ 0.001^{***} $  & $ 0.001^{***} $  & $ 0.001^{***} $  & $ 0.001^{***} $  \tabularnewline
 \hline 
2 & GARCH(1,1)-$t_5$ & $ 0.3636 $  & $ 0.007^{***} $  & $ 0.0320^{**} $  & $ 0.1968 $  & $ 0.0480^{**} $  & $ 0.0979^{*} $  \tabularnewline
 \hline                               
3 & AR(1)-GARCH(1,1) & $ 0.001^{***} $  & $ 0.001^{***} $  & $ 0.001^{***} $  & $ 0.002^{***} $  & $ 0.001^{***} $  & $ 0.001^{***} $  \tabularnewline
 \hline                                                        
4 & AR(1)-GARCH(1,1)-$t_5$  & $ 0.5834 $  & $ 0.3526 $  & $ 0.1748 $  & $ 0.5524 $  & $ 0.2278 $  & $ 0.2498 $  \tabularnewline
 \hline                                                
\hline
\end{tabular}
}\label{tab:baipval}
\footnotesize
\renewcommand\baselineskip{11pt}
\textbf{Note:} $p$-values of test statistics applied to GARCH-type models for NYSE monthly equal-weighted returns 1926.1 - 1999.12.
\end{table}

\section{TESTING MULTIVARIATE DISTRIBUTIONS}
\setcounter{theorem}{0}\setcounter{equation}{0}
Multivariate nonlinear dynamic models are gaining a lot of interest in econometrics literature. In this section, we briefly discuss how the results in this paper can be extended to the multivariate case. 
Suppose a sequence of $d\times 1$
vectors $Y_{1},Y_{2},\ldots,Y_{n}$, where $Y_{t}=(Y_{t1},Y_{t2},%
\ldots,Y_{td})^{\prime }, t=1,...,n,$ is given. Let $\Omega _{t}$ be again the information set at
time $t$, i.e. the $\sigma$-field of $%
\{Y_{t-1},Y_{t-2},\ldots \}$.
We are interested in the joint distributions $F_{t}(\cdot|\Omega _{t},\theta )$,
conditional on the past information, parameterized by $\theta \in \Theta $. The null hypothesis then is \bigskip

$H^{M}_{0}$ : The multivariate distribution of $Y_{t}$ conditional on $\Omega
_{t} $ is in the parametric family $F_{t}(\cdot\mid\Omega _{t},\theta )$ for some $%
\theta _{0}\in \Theta $.

For example, multivariate GARCH models (Engle, 2002), and dynamic copula models with parametric marginals (Patton, 2006), specify conditional joint distribution, and their specification can be tested using $H^{M}_{0}$. Following Rosenblatt (1952) and Diebold et al. (1999), define the multivariate probability integral transforms  as a univariate sequence  $U_{(t-1)d+k}=F_{tk}\left(Y_{tk}\mid\{Y_{t,k-1},\ldots,Y_{t,1},\Omega _{t}\},\theta_0\right)$, for $k=1,\ldots,d$, where $F_{tk}\left(\cdot\mid\{Y_{t,k-1},\ldots,Y_{t,1},\Omega _{t}\},\theta_0\right)$ is the distribution of $Y_{tk}$ conditional on $\{Y_{t,k-1},\ldots,Y_{t,1},\Omega _{t}\}$. This distribution can be computed from the null joint distribution $F_{t}(\cdot\mid\Omega _{t},\theta )$. Explicit formulas of the multivariate probability integral transforms are available for VAR and multivariate GARCH models with normal and student-$t$ innovations (Bai and Chen, 2009) and for dynamic copula models (Patton, 2013). 

Under $H^{M}_{0}$, the transforms $U_{\tau}$, which now constitute a univariate series of length $nd$, are iid uniform on $[0,1]$.  Thus we can apply test statistics based on $V_{2n,j}$ to this series, using that $P(U_{\tau}\le r_1,U_{\tau-j}\le r_2)=r_1 r_2$ for $j=1,2,\ldots$ under $H^{M}_{0}$. Note that for multivariate models, it is important to use multivariate process based tests even when the original multivariate series $Y_t$ is independent (across time) or iid. For instance, consider a bivariate normal iid series $Y_t$ with nonzero correlation between $Y_{t,1}$ and $Y_{t,2}$. If we test the null hypothesis that a series is iid normal and uncorrelated across two dimensions, then it will be hard to reject this null using only tests based on univariate distribution by the same reasoning as we discussed in the example of testing AR(1) against AR(2). Misspecification in cross sectional correlations is very undesirable. For instance, in financial applications such misspecification may result in overestimation of the effect of diversification strategies.

\section{CONCLUSION}
\setcounter{theorem}{0}\setcounter{equation}{0}
We fill a gap in the literature by introducing a test which inherits all the helpful features of tests based on the empirical process of generalized residuals and is consistent in a time series setup.
In some particular cases, it might be possible to further increase the power of the testing procedure. One possibility is to use a different functional $\Gamma$ and add weights to the process. 
The weak convergence result of the paper might be useful in these cases and could also assist in developing a theory for more general transformations, including using the martingale approach to get a distribution-free test, and multivariate models.

\section*{Acknowledgment}

I am grateful to Carlos Velasco for stimulating questions and  important suggestions and to Vanessa Berenguer, Victor Chernozhukov,  Miguel Delgado,  Manuel Dominguez, Jesus Gonzalo, Oliver Linton, Andrew Patton, Stefan Sperlich, Abderrahim Taamouti and the anonymous referees for helpful comments.  I thank Don Andrews, Yuichi Kitamura and Oliver Linton for possibility of visiting Cowles Foundation at Yale University and London School of Economics and their hospitality. I acknowledge financial support from the Spanish Ministerio de Economia y Competividad, Ref. no. ECO2012-31758.

\bibliographystyle{chicago}

\section*{Appendix A: General weak convergence result}
\renewcommand{\theequation}{A.\arabic{equation}}
\renewcommand{\thesection}{A}
\setcounter{equation}{0}
%
We introduce some notation and prove a
general weak convergence result for our process. For data generated under $%
G_{nt}$ with a true parameter denoted by $\theta _{n}$, $\hat{U}_{t}=F_{t}(Y_{t}|\hat{\theta})$ are not uniform iid (unless
$\delta =0$ and $\hat{\theta}=\theta _{n}$), but instead $U_{t}^{\dagger
}:=G_{nt}(Y_{t}|\theta _{n})$ are. So we have%
\begin{eqnarray*}
\hat{U}_{t}\leq r_{i} &\Leftrightarrow &F_{t}(Y_{t}|\hat{\theta})\leq r_{i}
\Leftrightarrow F_{t}(G_{nt}^{-1}(U_{t}^{\dagger }|\theta _{n})|\hat{\theta}%
)\leq r_{i} \\
&\Leftrightarrow &U_{t}^{\dagger }\leq G_{nt}(F_{t}^{-1}(r_{i}|\hat{\theta}%
)|\theta _{n}),
\end{eqnarray*}%
hence $I(\hat{U}_{t}\leq r_{i})=I(U_{t}^{\dagger }\leq \eta _{t}\left(
r_{i}\right) )$ where $\eta _{t}\left( r_{i}\right)=\eta_{t}\left( r_{i},\hat\theta\right)
=G_{nt}(F_{t}^{-1}(r_{i}|\hat{\theta})|\theta _{n})$. Define further $%
\Delta _{t}\left( a\right) $ $=I(U_{t}^{\dagger }\leq a).$ Then
\begin{equation*}
V_{2n}(r)=\frac{1}{\sqrt{n-1}}\sum_{t=2}^{n}\left[ \Delta _{t}\left(
r_{1}\right) \Delta _{t-1}\left( r_{2}\right) -r_{1}r_{2}\right]
\end{equation*}%
and%
\begin{equation*}
\hat{V}_{2n}(r)=\frac{1}{\sqrt{n-1}}\sum_{t=2}^{n}\left[ \Delta _{t}\left(
\eta _{t}\left( r_{1}\right) \right) \Delta _{t-1}\left( \eta _{t-1}\left(
r_{2}\right) \right) -r_{1}r_{2}\right] .
\end{equation*}%

We need the following conditions on $\eta _{t}.$ As it is discussed in the proofs of Propositions~3.3-3.6 below, these assumptions  hold for each $\eta _{t}$ considered in the paper and therefore  impose no additional restrictions to those listed in the main text of the paper. 
\begin{itemize}
\item[(C1)] 
\begin{equation*}
\sup_{r\in [0,1], t=1,..,n}\left\vert \eta _{t}\left( r\right) -r\right\vert
=O_{p}\left( n^{-1/2}\right).
\end{equation*}
\item[(C2)]  $\forall M\in(0,\infty)$, $\forall M_2\in(0,\infty)$ and $\forall\delta> 0 $
\begin{equation*}
\sup_{r\in [0,1]}
\frac{1}{\sqrt{n}}
\sum_{t=1}^{n}
\sup_{
\substack{
||u-v||\le
 M_2 n^{-1/2-\delta}\\
u,v\in B\left(\theta_n, M n^{-1/2}\right)
}
 }
\left|
{\eta
}_{t}\left( r,u\right)
-
{\eta }_{t}\left(
r,v\right)
\right| 
=o_{p}\left(1\right).
\end{equation*}
\item[(C3)]  $\forall M\in(0,\infty)$, $\forall M_2\in(0,\infty)$ and $\forall\delta> 0 $
\begin{equation*}
\sup_{|r-s|\le M_2 n^{-1/2-\delta} }
\frac{1}{\sqrt{n}}
\sum_{t=1}^{n}
\sup_{u\in B\left(\theta_n, M n^{-1/2}\right)}
\left|
{\eta
}_{t}\left(r, u\right)
-
{\eta }_{t}\left(s, u\right)
\right| 
=o_{p}\left(1\right).
\end{equation*}
\end{itemize}

\begin{lem}\label{lemmaz}
Let
\begin{equation*}
\xi_{t}= 
\left(
\Delta _{t-1}\left( \eta_{t-1}\left( r_{2}\right) \right)
-\Delta _{t-1}\left( r_{2}\right)  
\right)
\left( \Delta
_{t}\left( r_{1}\right) -r_{1}
\right)
\end{equation*}
and $z_{n}=\frac{1}{\sqrt{n}}\sum_{t=1}^{n}\xi_{t},$ then under (C1)-(C3), $\sup_{r}\left\vert z_{n}\right\vert =o_{p}(1).$
\end{lem}

\proofs{Lemma \ref{lemmaz}}{
Note that for any function $z_{n}(r,\hat{\theta}),$ depending on random
vector $\hat{\theta},$ with $\sqrt{n}\left( \hat{\theta}-\theta _{n}\right)
=o_{p}(1),$ for $\sup_{r}\left\vert z_{n}(r,\hat{\theta})\right\vert
=o_{p}(1)$ it is sufficient to show for some $\gamma <1/2,$
\begin{equation}\label{eq:ngamma}
\sup_{r,\left\Vert \eta -\theta _{n}\right\Vert \leq n^{-\gamma }}\left\vert
z_{n}(r,\eta )\right\vert =o_{p}(1).
\end{equation}%
Indeed,%
\begin{eqnarray*}
P\left( \sup_{r}\left\vert z_{n}(r,\hat{\theta})\right\vert >\varepsilon
\right) &\leq &P\left( \sup_{r,\left\Vert \eta -\theta _{n}\right\Vert \leq
n^{-\gamma }}\left\vert z_{n}(r,\eta )\right\vert >\varepsilon \right) \\
&&+P\left( \sqrt{n}\left\vert \hat{\theta}-\theta _{n}\right\vert
>n^{1/2-\gamma }\right),
\end{eqnarray*}%
where the second summand is $o\left( 1\right) $. From now on, all $\sup $
are taken with respect to $r$ and nonrandom $\eta $ s.t. $\left\Vert \eta
-\theta _{n}\right\Vert \leq n^{-\gamma }$.
We bound the expectation of the supremum with the
expectation of the maximum over a finite number of points, which itself is
bounded by the sum of the expectations. Having expectations ``inside" allow
to go from indicators to smooth functions, the difference of which can be bounded
by the differences of their arguments (see, e.g., Boldin, 1989).

First of all, we will show that $\forall\ \eta, r\ \left|z_{n}\right|=o_p\left(n^{-1/2}\right)$.
Since $\xi_{t}$  are bounded by 1 in absolute value and form a martingale
difference sequence with respect to $\Omega _{t}$,
 by the Doob
inequality  $\forall p\ge 1$ and $\forall \varepsilon>0$%
\begin{equation*}
P\left( \max_{t=1,
\ldots,n}\left\vert z_{t} \right\vert
>\varepsilon \right) \leq 
E \left| z_{n}\right| ^{p}/\varepsilon ^{p} .
\end{equation*}%
Next, by the Rosenthal 
inequality (Hall and Heyde, 1980, page 23),
$\forall p\ge 2\ \exists C_1$%
\begin{equation*}
E \left| z_{n}\right| ^{p} \leq n^{-p/2}C_1\left[
E\left\{ \sum E\left( \xi_t^{2}|\Omega
_{t}\right) \right\} ^{p/2}
+\sum E\left| \xi_t
\right| ^{p}\right].
\end{equation*}%
The conditional variance is given by 
\begin{equation*}
E\left(\xi_{t}^{2}|\Omega _{t}\right)
=
\left|\Delta _{t-1}\left({\eta }_{t-1}\left(r_{2}\right) \right)
-\Delta _{t-1}\left( r_{2}\right) \right|\left(r_1-r_1^2\right).
\end{equation*}%
By (C1) and since $\forall I\subset \{1,\ldots,n\}$, for $t'=\max I$ and
$\forall k_t \in \{1,\ldots,p/2\}$
\begin{eqnarray*}
&&E
\prod_{t\in I}\left|\Delta _{t}\left({\eta }_{t}\left(r_{2}\right) \right)
-\Delta _{t}\left( r_{2}\right) \right|^{k_t}
=E
\prod_{t\in I}
\left|\Delta _{t}\left({\eta }_{t}\left(r_{2}\right) \right)
-\Delta _{t}\left( r_{2}\right) \right|\\
&=&
E\left[\left\{
\prod_{t\in I\setminus\{t'\}} 
\left|\Delta _{t}\left({\eta }_{t}\left(r_{2}\right) \right)
-\Delta _{t}\left( r_{2}\right) \right| \right\}
E\left\{
\left|\Delta _{t'}\left({\eta }_{t'}\left(r_{2}\right) \right)
-\Delta _{t'}\left( r_{2}\right) \right| 
\mid\Omega_{t'}\right\}\right], \\
\end{eqnarray*}
which by uniformity and independence of $U_t^{\dag}$ of the past equals to
\begin{equation*}
E\left[\left\{
\prod_{t\in I\setminus\{t'\}} 
\left|\Delta _{t}\left({\eta }_{t}\left(r_{2}\right) \right)
-\Delta _{t}\left( r_{2}\right) \right| \right\}
\left|{\eta }_{t'}\left( r_{2}\right) - r_{2}\right|
\right],
\end{equation*}
which by Holder inequality for $1<q<+\infty$ is bounded by
\begin{equation*}
\left[E\left\{
\prod_{t\in I\setminus\{t'\}} 
\left|\Delta _{t}\left({\eta }_{t}\left(r_{2}\right) \right)
-\Delta _{t}\left( r_{2}\right) \right| \right\}^{q}\right]^{\frac{1}{q}}
\left[E\left|{\eta }_{t'}\left( r_{2}\right) - r_{2}\right|^{\frac{q}{q-1}}\right]^{1-\frac{1}{q}}.
\end{equation*}
The second term is $O\left(n^{-1/2}\right)$. We repeat this inequality for 
$I_1=I\setminus\{t'\}$ and so on. Using multinomial formula and $q$ sufficiently close to $1$, and taking $p=6$, we get
\begin{equation*}
n^{-p/2}E\left\{ \sum E\left( \xi_{t}^{2}|\Omega
_{t}\right) \right\} ^{p/2}=O\left(n^{-3/2}\right).
\end{equation*}
Because of boundedness of the indicator, $\left|\xi_{t}\right|\le 1$ and
\begin{equation*}
n^{-p/2}\sum E\left|\xi_{t}\right|^{p}
=O\left(n^{1-p/2}\right).
\end{equation*}
Thus,
\begin{equation}\label{eq:nuniformmds}
P\left( \max_{t=1,
\ldots,n}\left\vert z_{t} \right\vert
>\varepsilon \right) = O\left(n^{-3/2}\right).
\end{equation}

Now we establish the uniform result. Break up the interval 
$[-n^{-\gamma },n^{-\gamma }]$ ($\gamma$  defined in Equation \ref{eq:ngamma}) into ${m_{n}}$ parts
by the points $-n^{-\gamma }+2n^{-\gamma }s/{m_{n}},\ s=1,...,{m_{n}}.$
These points split the cube $[-n^{-\gamma },n^{-\gamma }]^{q}$ into $%
{m_{n}^q}$ cubes with vertexes at these points. Now, at the intersection of
each cube $s=1,...,{m_{n}q}$ with the sphere $\left\Vert \eta _{n}-\theta
_{n}\right\Vert $ $\leq n^{-\gamma }$, denote maximum and minimum of $\eta
_{t}\left( r_{i}\right) $ as $\overline{\eta }_{ts}\left( r_{i}\right) $ and
$\underline{\eta }_{ts}\left( r_{i}\right) $. Divide also the interval $%
[0,1] $ into $N_{n}$ equal intervals. Let $\left[
r_{i}^{s_{i}-1},r_{i}^{s_{i}}\right] ,s_{i}=1,...,N_{n,}$ denote the
interval which contains point $r_{i}$.
Since the indicator is monotonous, if $\eta _{n}$ is in the cube $s$, then $\sup
z_{n}$ is bounded from above by the maximum over $s,s_{1},s_{2}$ of
\begin{eqnarray}
&&\frac{1}{\sqrt{n}}\sum_{t=2}^{n}
\left\{ 
\Delta _{t-1}\left( \overline{\eta }%
_{t-1s}\left( r_{2}^{s_{2}}\right) \right) 
 -\Delta _{t-1}\left(
r_{2}^{s_{2}}\right) 
\right\}
\left\{
\Delta _{t}\left( r_{1}^{s_{1}}\right)
-r_{1}^{s_{1}}
\right\}
\label{eq:discr-mds} \\
&&+\frac{1}{\sqrt{n}}\sum_{t=2}^{n}\left\{ \Delta _{t-1}\left(
r_{2}^{s_{2}}\right) 
 \Delta _{t}\left( r_{1}^{s_{1}}\right)
 -\Delta _{t-1}\left( r_{2}^{s_{2}-1}\right) 
\Delta _{t}\left( r_{1}^{s_{1}-1}\right) \right\}  \label{eq:discr-proc} \\
&&+\frac{1}{\sqrt{n}}\sum_{t=2}^{n}\left\{ \Delta _{t-1}\left( \overline{\eta
}_{t-1s}\left( r_{2}^{s_{2}}\right) \right) 
r_{1}^{s_{1}}
-
\Delta _{t-1}\left( \underline{\eta }_{t-1s}\left(
r_{2}^{s_{2}-1}\right) \right) 
r_{1}^{s_{1}-1}\right\}. \label{eq:discr-eta} 
\end{eqnarray}%
The maximum of the absolute value of Equation (\ref{eq:discr-mds}) is denoted  by $z_{n}\left( s,s_{1},s_{2}\right) $ and treated using (\ref{eq:nuniformmds})
\begin{eqnarray*}
P\left( \max_{s,s_{1},s_{2}}\left\vert z_{n}\left( s,s_{1},s_{2}\right)
\right\vert >\varepsilon \right) 
&\leq &
\sum_{s,s_{1},s_{2}}P\left(
\left\vert z_{n}\left( s,s_{1},s_{2}\right) \right\vert >\varepsilon \right)
\\
&\leq &\sum_{s,s_{1},s_{2}}{E\left( z_{n}\left( s,s_{1},s_{2}\right)
^{2}\right) }/{\varepsilon ^{2}} \\
&=&O\left( {{m_{n}^q}N_{n}^{2}}{n^{-3/2}}\right) .
\end{eqnarray*}%
The maximum of the absolute value of Equation (\ref{eq:discr-proc}), is no more  than
\begin{eqnarray*}
\sup_{\left|r_i-s_i\right|\le N_n^{-1}}\left|
V_{2n}{\left(r\right)}
-
V_{2n}{\left(s\right)}
\right| =O_p(N_n^{-1}).
\end{eqnarray*}%
Finally, Equation (\ref{eq:discr-eta}) equals
\begin{eqnarray*}
&&\frac{1}{\sqrt{n}}\sum_{t=2}^{n}
r_{1}^{s_{1}}
\left\{ \Delta _{t-1}\left( \overline{\eta
}_{t-1s}\left( r_{2}^{s_{2}}\right) \right) 
-
\Delta _{t-1}\left( \underline{\eta }_{t-1s}\left(
r_{2}^{s_{2}-1}\right) \right) 
-
\overline{\eta
}_{t-1s}\left( r_{2}^{s_{2}}\right)
+
\underline{\eta }_{t-1s}\left(
r_{2}^{s_{2}-1}\right)
\right\} \\
&&+\frac{1}{\sqrt{n}}\sum_{t=2}^{n}
r_{1}^{s_{1}}
\left\{ 
\overline{\eta
}_{t-1s}\left( r_{2}^{s_{2}}\right)
-
\underline{\eta }_{t-1s}\left(
r_{2}^{s_{2}-1}\right)
\right\} \\
&&+\frac{1}{\sqrt{n}}\sum_{t=2}^{n} 
\Delta _{t-1}\left( \underline{\eta }_{t-1s}\left(
r_{2}^{s_{2}-1}\right) \right) 
\left\{r_{1}^{s_{1}}-r_{1}^{s_{1}-1}\right\} .
\end{eqnarray*}%
The maximum of the absolute value of the first term is less than 
\begin{equation*}
\sup_{|r-s|\to 0 }\left|
\frac{r_{1}^{s_{1}}}{\sqrt{n}}
\sum_{t=1}^{n}
\left\{ \Delta _{t}\left( r \right) 
-
\Delta _{t}\left(s \right) 
-r+s
\right\}
\right|=o_p(1).
\end{equation*}
The second term equals 
\begin{equation*}
\frac{r_{1}^{s_{1}}}{\sqrt{n}}
\sum_{t=1}^{n}
\left\{ 
\overline{\eta
}_{t-1s}\left( r_{2}^{s_{2}-1}\right)
-
\underline{\eta }_{t-1s}\left(
r_{2}^{s_{2}-1}\right)
\right\}
+
\frac{r_{1}^{s_{1}}}{\sqrt{n}}
\sum_{t=1}^{n}
\left\{ 
\overline{\eta
}_{t-1s}\left( r_{2}^{s_{2}}\right)
-
\overline{\eta }_{t-1s}\left(
r_{2}^{s_{2}-1}\right)
\right\}
\end{equation*}
and the maximum of its absolute value is 
$O_p(n^{-\gamma+1/2}m_n^{-1}\sqrt{q})+O_p\left(N_n^{-1} n^{1/2}\right)$ %
 by (C2) and (C3). The third term is $O_p\left(N_n^{-1} n^{1/2}\right)$.
If we take $N_n^2=n^{1+\delta_N}$ and ${m_n^q}=n^{1/2-\delta_m}$, where
$\delta_N=\delta_m /2=\left(\gamma q+1/2-q/2\right)/2$,
all parts will be $o_p(1)$ for $1/2-1/(2q)<\gamma<1/2$.
The same argument holds for the bound of $\sup z_n$ from below.  } 

Define
\begin{equation*}
k_{n}=\frac{1}{\sqrt{n}}\sum_{t=2}^{n}\left[ \Delta _{t-1}\left(
\eta _{t-1}\left( r_{2}\right)\right) \left( \eta _{t}\left( r_{1}\right) -r_{1}\right) +r_{1}\left(
\eta _{t-1}\left( r_{2}\right) -r_{2}\right) \right],
\end{equation*}%

Now we are able to state a general result on the parameter estimation effect on the asymptotics of our process:

\begin{lem}
\label{lem:weakconv} Under (C1)-(C3)
\begin{equation}\label{eq:implemma}
\sup_{r}\left\vert \hat{V}_{2n}(r)-V_{2n}(r)-k_{n}\right\vert =o_{p}(1).
\end{equation}
\end{lem}

\proofs{Lemma \ref{lem:weakconv}}{
Follows from Lemma \ref{lemmaz} since
\begin{eqnarray*}
&&\hat{V}_{2n}(r)-V_{2n}(r)-k_{n}-z_n\\
&=&\frac{1}{\sqrt{n}}\sum_{t=2}^{n}
\Delta _{t-1}\left( \eta_{t-1}\left( r_{2}\right) \right) 
\left\{ 
\Delta _{t}\left( 
\eta_{t}\left( r_1\right) 
\right) 
-
\Delta _{t}\left( 
r_{1} \right) 
-
\eta_{t}\left( r_{1}\right)
+
r_1
\right\}  \\
&&+\frac{1}{\sqrt{n}}\sum_{t=2}^{n}
r_1
\left\{ 
\Delta _{t-1}\left( 
\eta_{t-1}\left( r_2\right) 
\right) 
-
\Delta _{t-1}\left( 
r_{2} \right) 
-
\eta_{t-1}\left( r_{2}\right)
+
r_2
\right\}  ,
\end{eqnarray*}%
and the maximum of its absolute value is $o_p(1)$, by the argument similar to Koul (1996).
}%

\section*{Appendix B: Proofs of Propositions}
\renewcommand{\theequation}{B.\arabic{equation}}
\renewcommand{\thesection}{B}
\setcounter{equation}{0}

\proofs{Proposition \ref{propFY}}{The proof is standard and thus omitted.}

\proofs{Proposition \ref{proplimitV}}
{We use the functional CLT of Pollard (1984, Theorem 10.12). We need to check
equicontinuity and convergence of finite dimensional distributions. Equicontinuity can be shown in a standard way. 

We verify that the process has zero mean and covariance converging to (\ref%
{eqasscov}). The finite dimensional distributions converge by the Cramer-Wold device and a CLT for stationary 1-dependent data with a finite third moment (see Theorem 2 in Hoeffding and Robbins, 1948). Indeed, because of the independence and uniformity  of $U_{t}$,  the  mean is
\begin{equation*}
EV_{2n}(r) =\frac{1}{\sqrt{n-1}}\sum_{t=2}^{n}\left[ P(U_{t}\leq r_{1},U_{t-1}\leq
r_{2})-r_{1}r_{2}\right]  =0,
\end{equation*}
and variance is derived using
\begin{eqnarray*}
&&E\left[ I(U_{t}\leq r_{1},U_{t-1}\leq r_{2})-r_{1}r_{2}\right] \left[
I(U_{t^{\prime }}\leq s_{1},U_{t^{\prime }-1}\leq s_{2})-s_{1}s_{2}\right] \\
&=&P(U_{t}\leq r_{1},U_{t-1}\leq r_{2},U_{t^{\prime }}\leq
s_{1},U_{t^{\prime }-1}\leq s_{2})-r_{1}r_{2}s_{1}s_{2} \\
&=&\left\{
\begin{array}{cc}
0, & \text{if }|t-t^{\prime }|>1
\\
(r_{1}\wedge s_{1})(r_{2}\wedge s_{2})-r_{1}r_{2}s_{1}s_{2}, & \text{%
if }t^{\prime }=t \\
(r_{1}\wedge s_{2})r_{2}s_{1}-r_{1}r_{2}s_{1}s_{2}, & \text{if }%
t^{\prime }=t+1 \\
(r_{2}\wedge s_{1})r_{1}s_{2}-r_{1}r_{2}s_{1}s_{2}, & \text{if }%
t^{\prime }=t-1.%
\end{array}%
\right.
\end{eqnarray*}%
Then covariance of the process is%
\begin{align*} 
E\left[V_{2n}(r)V_{2n}(s)\right]
\overset{n\rightarrow \infty }{\rightarrow }&
(r_{1}\wedge s_{1})(r_{2}\wedge s_{2})-r_{1}r_{2}s_{1}s_{2}\\
&+
(r_{1}\wedge s_{2})r_{2}s_{1}-r_{1}r_{2}s_{1}s_{2}+
(r_{2}\wedge s_{1})r_{1}s_{2}-r_{1}r_{2}s_{1}s_{2}.
\end{align*}
}

Before we move on, we derive the joint convergence of ${V_{2n}(r)}$ and ${\sqrt{n}(\hat{%
\theta}}-{\theta }_{n}).$

\begin{lem}
\label{lemmaVpsiasB} Under Assumptions 3.1-3.5 (%
3.1, 3.2, 3.4-3.6) under $G_{nt}$ (under $\{\theta _{n}:n\geq 1\}$) we have%
\begin{equation*}
\binom{{V_{2n}}}{{\sqrt{n}(\hat{\theta}}-{\theta }_{n})}
\weakc
\binom{{V_{2\infty }}}{\psi _{\infty }+\mu}.
\end{equation*}
\end{lem}

\proofs{Lemma \ref{lemmaVpsiasB}}{
Assume $\mu=0$, otherwise subtract it from the left hand side.  By Proposition \ref{proplimitV} and CLT for MDS from the expansion
of Assumption 3.3,
we have componentwise convergence. 
To prove vector convergence we use functional CLT of Pollard (1984, Theorem
10.12).
We need to check equicontinuity and convergence of finite dimensional distributions.
Equicontinuity follows from the fact that equicontinuity of the vector is
equivalent to equicontinuity of its components and the equicontinuity of the first
component provided in Proposition \ref{proplimitV}. Other components are equicontinuous automatically since they do not depend on
parameter~$r$.

To check the convergence of finite dimensional distributions, we apply the Cramer-Wold device 
and prove a CLT for triangular arrays $\sum^n_{t=1}  v_t$, where
\begin{equation*}
v_t=\frac{1}{
\sqrt{n-1}}\lambda'
\begin{pmatrix}
{I({U}_{nt}\leq r^1_{1})I({U}_{nt-1}\leq r^1_{2})-r^1_{1}r^1_{2}} \\
\vdots\\
{I({U}_{nt}\leq r^{k}_{1})I({U}_{nt-1}\leq r^{k}_{2})-r^{k}_{1}r^{k}_{2}}\\
{\psi\left(Y_{nt},\Omega_{nt},\theta_n\right)}
\end{pmatrix},
\end{equation*}
for any $k$, a column vector $\lambda\in R^{k+L}$, s.t. $\lambda'\lambda=1$, and any $r^j_i\in [0,1]$ with $j=1,\ldots,k$ and $i=1,2$. 
We derive a CLT for sum of $v_{t}$, whose summands $1,\ldots, k$ form $1$-dependent sequences,
while summands $k+1,\ldots, k+L$ form MDS.
Also note, that sequence $v_t$ is lag-$1$ serially correlated.

To prove the  CLT,   we split $v_t$ into groups, skipping $1$ element to
form  a sum of MDS. Fix $p>1$ and denote $m=\lfloor \frac{n}{p}\rfloor
$. Consider the following decomposition%
\begin{equation*}
\sum^n_{t=1} v_t=A_{n,p}+B_{n,p}+C_{n,p},
\end{equation*}
where $A_{n,p}=\sum_{k=1}^{m}A_{n,p,k}$ denotes the sum
of blocks $A_{n,p,k}=\sum_{t=(k-1)p+2}^{kp}v_{t}$ of length $p-1$, $B_{n,p}=%
\sum_{i=1}^{\left( m-1\right) }v_{ip+1}$ and 
$C_{n,p}=\sum_{t=mp+1}^{n}v_{t}$ 
denotes the sum of remaining $n-pm<p$ terms. We will show
now that $A_{n,p}$ converges to the right limit, $B_{n,p}$ and $C_{n,p}$ are
$o_{p}(1)$. 

Considering that $B_{n,p}$ and $C_{n,p}$ are uncorrelated, 
$B_{n,p}$ are serially uncorrelated and $C_{n,p}$ has no more than $p-1$ $1$-dependent terms, we have
\begin{eqnarray*}
P\left( |B_{n,p}+C_{n,p}|>\varepsilon \right) &\leq &
\frac{1}{\varepsilon ^{2}} E\left( B_{n,p}+C_{n,p}\right) ^{2} = \frac{1}{\varepsilon ^{2}}\left[ EB_{n,p}^{2}+EC_{n,p}^{2}%
\right]  \\
&= &\frac{1}{\varepsilon ^{2}}\left[ \sum_{i=1}^{m-1}
Ev_{ip+1}^{2 } \right] +
\frac{1}{\varepsilon ^{2}}
\left [\sum_{t=mp+1}^{n} Ev_{t}^{2}
+2\sum_{t=mp+2}^{n} Ev_{t}v_{t-1}\right].
\end{eqnarray*}%
Allowing $m\to\infty$ and $p\to\infty$,  
and since $n$ grows faster than $m$ and $p$ and  
the summands are uniformly bounded by Assumption 3.3 
and Cauchy-Schwarz inequality, we can always make this quantity 
arbitrarily small.

$A_{n,p,k}$ is an MDS with respect
to the new filtration $\Omega _{k}^{\prime }=\Omega _{(k-1)p+1}$. 
To show the asymptotic normality of  $A_{n,p}$
we apply CLT for MDS (Hall and Hyde, 1980, Corollary 3.1), 
which requires two conditions:
1) Lindeberg condition and 
2) convergence of conditional variances. They can be verified for each summand
and follow from the boundedness of indicators, 
the conditions in Assumption 3.3 and 
the convergence of the unconditional variance.
 The variance is approximated by
\begin{equation*}
EA_{n,p,k}^{2} =\sum_{k=1}^{m}\left(
\sum_{t=(k-1)p+2}^{kp}Ev_{t}^{2 }+2%
\sum_{t=(k-1)p+3}^{kp}Ev_{t}v_{t-1}\right),
\end{equation*}%
which coincides in the limit with
\begin{equation*}
E\left(\sum_{k=1}^n v_t\right)^{2} =\sum_{k=1}^{n}Ev_{t}^{2 }+2%
\sum_{t=2}^{n}Ev_{t}v_{t-1}.
\end{equation*}%
Thus we have shown the convergence of the finite dimensional distributions.

Now we write the asymptotic covariance matrix.
The covariance of the process $V_{2\infty}$ is derived in Proposition \ref{proplimitV}, 
see Equation (\ref{eqasscov}). 
The other component on the main diagonal is $\Psi$ as in Assumption 3.3. 
The covariance between the two components, under $\{\theta
_{n}:n\geq 1\}$, is
\begin{eqnarray*}
&&\cov(V_{2n}(r),\frac{1}{\sqrt{n}}\sum_{t=2}^{n}\psi\left(Y_{nt},\Omega_{nt},\theta_n\right)) \\
&&=\frac{1}{n}\sum_{t=2}^{n}E\left[ \left( I(U_{nt}\leq r_{1})I(U_{nt-1}\leq
r_{2})+r_{1}I(U_{nt}\leq r_{2})\right) \psi\left(Y_{nt},\Omega_{nt},\theta_n\right)\right].
\end{eqnarray*}
}

\proofs{Proposition \ref{propnull}} {
We use Lemmas \ref%
{lem:weakconv} and \ref{lemmaVpsiasB} with $\theta _{n}=\theta _{0}$ and $\delta =0,$ i.e. with $%
\eta _{t}\left( r_{i}\right) =F_{t}(F_{t}^{-1}(r_{i}|\hat{\theta})|\theta
_{0}).$ Conditions (C1)-(C3) follow from Assumption 3.2. Condition 
(d) ensures $\sup_{r}\left\vert k_{n}-h(r)\right\vert
=o_{p}(1).$ The result follows from the functional continuous mapping theorem (CMT) (Pollard, 1984, Theorem IV.12, p.70).}%

\proofs{Proposition \ref{propalt}}{
We use Lemmas \ref%
{lem:weakconv} and \ref{lemmaVpsiasB}  with $\theta _{n}=\theta _{0}$ and $\delta \not=0,$ i.e. with
$\eta _{t}\left( r_{i}\right) =F_{t}(F_{t}^{-1}(r_{i}|\hat{\theta})|\theta
_{0})+\left( H_{t}(F_{t}^{-1}(r_{i}|\hat{\theta}%
))-F_{t}(F_{t}^{-1}(r_{i}|\hat{\theta})|\theta _{0})\right){\delta }/{\sqrt{n}}.$
Conditions (C1)-(C3) follow from Assumption 3.2 and continuity of $H_t$, Assumption 3.4, 
since $\left\vert \eta
_{t}\left( r_{i}\right) -F_{t}(F_{t}^{-1}(r_{i}|\hat{\theta})|\theta
_{0})\right\vert \leq {\delta }/{\sqrt{n}}$. Condition (d) 
ensures that $\sup_{r}\left\vert k_{n}-h(r)\right\vert =o_{p}(1).$ The result
follows from CMT.}

\proofs{Proposition \ref{propFA}}{
We use Lemma \ref%
{lem:weakconv}   with  $\theta _{n}= \theta _{1}$ and $%
\delta =0,$ i.e. with $\eta _{t}\left( r_{i}\right) =H_{t}(F_{t}^{-1}(r_{i}|%
\hat{\theta})|\theta _{1}).$ Conditions (C1)-(C3) follow from Assumption 3.2. Assumption 3.6 controls the drift from the parameter estimation.  The result follows from CMT.}

\proofs{Proposition \ref{propboot}}{
We use Lemmas \ref%
{lem:weakconv} and \ref{lemmaVpsiasB}  with nonrandom $\theta _{n}\rightarrow \theta _{0}$ and $%
\delta =0,$ i.e. with $\eta _{t}\left( r_{i}\right) =F_{t}(F_{t}^{-1}(r_{i}|%
\hat{\theta})|\theta _{n}).$ Conditions (C1)-(C3) follow from Assumption 3.2. Condition (d)  ensures $%
\sup_{r}\left\vert k_{n}-h(r)\right\vert =o_{p}(1).$ The result follows from
CMT.}

\end{document}